\documentclass[12pt,epsf]{article}

\usepackage{times}
\usepackage{graphicx}
\usepackage{amsfonts}
\usepackage{amsmath}
\usepackage{amssymb}
\usepackage{amstext}
\usepackage{amsthm}
\usepackage{youngtab}

\usepackage{epsfig}
\usepackage{color}
\usepackage{leftidx}
\usepackage{tikz}
\psfigdriver{dvips}
\usepackage{latexsym}
\usepackage[matrix,curve]{xypic}
\usepackage{rotating}
\rotdriver{dvips}

\begin{document}

\baselineskip 6mm
\renewcommand{\thefootnote}{\fnsymbol{footnote}}


\newcommand{\nc}{\newcommand}
\newcommand{\rnc}{\renewcommand}


\rnc{\baselinestretch}{1.24}    
\setlength{\jot}{6pt}       
\rnc{\arraystretch}{1.24}   

\makeatletter
\rnc{\theequation}{\thesection.\arabic{equation}}
\@addtoreset{equation}{section}
\makeatother



\nc{\be}{\begin{equation}}

\nc{\ee}{\end{equation}}

\nc{\bea}{\begin{eqnarray}}

\nc{\eea}{\end{eqnarray}}

\nc{\xx}{\nonumber\\}

\nc{\ct}{\cite}

\nc{\la}{\label}

\nc{\eq}[1]{(\ref{#1})}

\nc{\newcaption}[1]{\centerline{\parbox{6in}{\caption{#1}}}}

\nc{\fig}[3]{
\begin{figure}
\centerline{\epsfxsize=#1\epsfbox{#2.eps}}
\newcaption{#3. \label{#2}}
\end{figure}
}


\def\CA{{\cal A}}
\def\CC{{\cal C}}
\def\CD{{\cal D}}
\def\CE{{\cal E}}
\def\CF{{\cal F}}
\def\CG{{\cal G}}
\def\CH{{\cal H}}
\def\CK{{\cal K}}
\def\CL{{\cal L}}
\def\CM{{\cal M}}
\def\CN{{\cal N}}
\def\CO{{\cal O}}
\def\CP{{\cal P}}
\def\CR{{\cal R}}
\def\CS{{\cal S}}
\def\CU{{\cal U}}
\def\CV{{\cal V}}
\def\CW{{\cal W}}
\def\CY{{\cal Y}}
\def\CZ{{\cal Z}}


\def\IB{{\hbox{{\rm I}\kern-.2em\hbox{\rm B}}}}
\def\IC{\,\,{\hbox{{\rm I}\kern-.50em\hbox{\bf C}}}}
\def\ID{{\hbox{{\rm I}\kern-.2em\hbox{\rm D}}}}
\def\IF{{\hbox{{\rm I}\kern-.2em\hbox{\rm F}}}}
\def\IH{{\hbox{{\rm I}\kern-.2em\hbox{\rm H}}}}
\def\IN{{\hbox{{\rm I}\kern-.2em\hbox{\rm N}}}}
\def\IP{{\hbox{{\rm I}\kern-.2em\hbox{\rm P}}}}
\def\IR{{\hbox{{\rm I}\kern-.2em\hbox{\rm R}}}}
\def\IZ{{\hbox{{\rm Z}\kern-.4em\hbox{\rm Z}}}}


\def\a{\alpha}
\def\b{\beta}
\def\d{\delta}
\def\ep{\epsilon}
\def\ga{\gamma}
\def\k{\kappa}
\def\l{\lambda}
\def\s{\sigma}
\def\t{\theta}
\def\w{\omega}
\def\G{\Gamma}


\def\half{\frac{1}{2}}
\def\dint#1#2{\int\limits_{#1}^{#2}}
\def\goto{\rightarrow}
\def\para{\parallel}
\def\brac#1{\langle #1 \rangle}
\def\curl{\nabla\times}
\def\div{\nabla\cdot}
\def\p{\partial}


\def\Tr{{\rm Tr}\,}
\def\det{{\rm det}}


\def\vare{\varepsilon}
\def\zbar{\bar{z}}
\def\wbar{\bar{w}}


\def\ad{\dot{a}}
\def\bd{\dot{b}}
\def\cd{\dot{c}}
\def\dd{\dot{d}}
\def\so{SO(4)}
\def\bfr{{\bf R}}
\def\bfc{{\bf C}}
\def\bfz{{\bf Z}}

\begin{titlepage}


\hfill\parbox{3.7cm} {{\tt arXiv:2309.05335}}

\vspace{15mm}

\begin{center}
{\Large \bf  Einstein Structure of Four-Manifolds}

\vspace{10mm}

Jeongwon Ho ${}^{a}$\footnote{freejwho@gmail.com}, Kyung Kiu Kim ${}^{b}$\footnote{kimkyungkiu@kookmin.ac.kr},
and Hyun Seok Yang ${}^c$\footnote{hsyang@gist.ac.kr}
\\[10mm]

${}^a$ {\sl Center for Quantum Spacetime, Sogang University, Seoul 121-741, Korea}

${}^b$ {\sl College of General Education, Kookmin University, Seoul 02707, Korea}

${}^c$ {\sl Department of Physics and Photon Science, Gwangju Institute of Science and Technology, Gwangju 61005, Korea}

\end{center}

\thispagestyle{empty}

\vskip1cm


\centerline{\bf ABSTRACT}
\vskip 4mm
\noindent

It is known that the moduli space of Einstein structures in four dimensions is generally considered to be rigid so that
Einstein metrics tend to be isolated modulo diffeomorphisms under infinitesimal Einstein deformations.
We examine the rigidity of the Einstein structure by considering deformations of the round four-sphere.
We show that any deviation from the standard metric of the round four-sphere (except for scaling) breaks
the Einstein condition. This further supports the idea of rigidity.
We analyze the Einstein structure of four-manifolds based on the irreducible decomposition
of the self-dual structure of Einstein manifolds.
\\


Keywords: Einstein Structure, Four-Manifold, Moduli Space

\vspace{1cm}

\today

\end{titlepage}

\renewcommand{\thefootnote}{\arabic{footnote}}
\setcounter{footnote}{0}

\section{Introduction}


An Einstein manifold has the property that its Ricci curvature tensor is proportional to its metric tensor.
Finding a metric that satisfies this condition while being nonhomogeneous and compact is quite challenging.
It requires careful balancing of the geometric properties to achieve the desired curvature.
Hence, it is difficult to find explicit solutions to Einstein's equations for compact nonhomogeneous manifolds
since the field equations for an appropriate metric ansatz are highly nonlinear partial differential equations.
For a compact Einstein manifold with positive Ricci scalar, the manifold must be simply connected.
For example, $\mathbb{S}^1 \times \mathbb{S}^3$ does not admit Einstein metrics \cite{hitchin1974}.
It turns out that there are not many simply connected compact manifolds.
Compact manifolds have a finite volume, which introduces constraints on their topology and geometry.
The moduli space of Einstein structures is, up to a discrete group action, a finite-dimensional manifold
(see Chap. 12 in \cite{besse}).
The rigidity of Einstein structures implies that under specific conditions, the only solutions to Einstein's equations
are homogeneous ones, leaving little room for nonhomogeneous constructions.
For example, G. Jensen proved that a homogeneous four-dimensional Einstein manifold is isometric to a Riemannian
symmetric space \cite{g.jensen}.
While some nonhomogeneous examples exist \cite{page,lebrun-jams}, they are generally not easy to find or generalize.
A nonhomogeneous space possesses less global symmetries, making it more difficult to find suitable geometric structures
that satisfy Einstein equations.

There are few known examples of compact four-dimensional symmetric Einstein manifolds
with positive scalar curvature: $\mathbb{S}^4, \, \mathbb{C}\mathbb{P}^2, \, \mathbb{S}^2 \times \mathbb{S}^2$.
These examples are all homogeneous manifolds. It was shown in \cite{hitchin1974} that
a connected sum such as $M = n \mathbb{C}\mathbb{P}^2$ allows an Einstein metric
only for $1 \leq n \leq 3$.\footnote{The explicit metrics
on $n \mathbb{C}\mathbb{P}^2$ in \cite{lebrun-jdg} are four-dimensional K\"ahler metrics with vanishing Ricci scalar,
which is not a vacuum solution but a solution in Einstein-Maxwell theory.}
There is also a compact nonhomogeneous Einstein manifold $(M, g)$ known as the Page space \cite{page}.
It is a solution of Einstein equations, $R_{ab} = \lambda \delta_{ab}$ with $\lambda >0$,
which is a nontrivial $\mathbb{S}^2$ bundle over $\mathbb{S}^2$.
It can be constructed as $M = \mathbb{C}\mathbb{P}^2 \sharp \, \overline{\mathbb{C}\mathbb{P}}^2$
by attaching the complex projective space and its complex conjugate,
removing a 4-ball from each manifold, and gluing them together
along the resulting $\mathbb{S}^3$ boundaries \cite{besse}.
This construction was generalized to $M = \mathbb{C}\mathbb{P}^2 \sharp \, n \overline{\mathbb{C}\mathbb{P}}^2$
for $1 \leq n \leq 8$ \cite{tian-yau,tian-im}. Their explicit metrics are not known in most cases
except for $n=1$ \cite{page}.
Only a finite number of compact four-manifolds may carry an Einstein metric,
which is conformal to an extremal K\"ahler metric \cite{besse,v21}.

In this paper, we will examine the rigidity of the Einstein structure by considering deformations
of the round sphere $\mathbb{S}^4$. The rigidity means that an Einstein structure has no infinitesimal Einstein deformations
and so it is an isolated point of the moduli space.
See $\mathbf{12.H}$ in \cite{besse} and the review \cite{rev-anderson} for the moduli space of Einstein structures.
The rigidity result for Einstein manifolds with positive Ricci scalar was shown in \cite{c-torre,torre-jmp} by proving that
there are no nontrivial solutions to the linearized instanton equations on conformally anti-self-dual Einstein spaces
with a positive cosmological constant.
It is also known \cite{besse,daga-yang} (see also \cite{v22,v23}) that the standard metrics on $\mathbb{S}^4$ or $\mathbb{C}\mathbb{P}^2$ are locally rigid,
i.e. the metrics are isolated points in the moduli space of Einstein metrics.
In order to study the moduli space of Einstein structures, we will use the irreducible decomposition of
curvature tensors \cite{ahs-sd4} (see also Chaps. $\mathbf{1.G}$ and $\mathbf{1.H}$ in Ref. \cite{besse}
and Secs. 1.1 and 2.1 in Ref. \cite{dona-kron}).
In section 2, we recapitulate the canonical decomposition of the curvature tensor in Refs. \cite{yang-col1,oh-yang,yang-col2,yang-col3,yang-col4} and provide a simple proof of
the Hitchin-Thorpe inequality \cite{hitchin1974,rev-anderson,pps-blms}.
We apply the canonical decomposition to the Page metric \cite{page} which is the first example of
a nonhomogeneous compact Einstein manifold and is a Hermitian-Einstein, conformal to a K\"ahler metric,
but is not K\"ahler itself \cite{besse,rev-anderson}.
In section 3, we examine the rigidity of the Einstein structure by considering deformations
of the round four-sphere $\mathbb{S}^4$.
We show that any deformation of the round four-sphere $\mathbb{S}^4$ causes it to deviate from the Einstein structure,
except trivial deformations at most only changing the size of the sphere.
This analysis may shed light on why Einstein metrics are typically isolated modulo diffeomorphisms,
implying a form of rigidity.
In section 4, we consider asymptotically locally flat (ALF) spaces which share the same asymptotic boundary conditions
but carry different topological structures.
Our formalism proves particularly useful for these spaces to reveal their Einstein structures.
In section 5, we derive an underlying equation governing the deformations of Einstein metrics that
may shed light on the physical origin of the rigidity of Einstein structures.
In appendix A, we will prove the strong version of the Hitchin-Thorpe inequality (Theorem $\mathbf{2}$ in \cite{hitchin1974})
in order to demonstrate the effectiveness of our approach.

\section{Einstein Structure of Four-Manifolds}

On an orientable Riemannian four-manifold, the vector space of 2-forms $\Omega^2 (M) = \Lambda^2 T^* M$ is decomposed into
the space of self-dual and anti-self-dual 2-forms \cite{dona-kron}
\begin{equation}\label{sd-asd}
  \Omega^2 (M) = \Omega_+^2 \oplus \Omega_-^2
\end{equation}
using the projection operators
\begin{equation}\label{proj-1}
  P_\pm = \frac{1}{2} (1 \pm *), \qquad P_\pm^2 = P_\pm, \qquad P_+ + P_- = 1,
\end{equation}
where $\Omega_\pm^2 \equiv  P_\pm \Omega^2$ is the $\pm 1$ eigenspaces of the Hodge star operator $*:  \Omega^2 \to \Omega^2$.
The Lorentz group $G = Spin(4)$ is a double cover of the four-dimensional Euclidean Lorentz group $SO(4)$,
 i.e., $SO(4) \cong Spin(4)/\mathbb{Z}_2$.
The Lorentz group $Spin(4)$ splits into a product of two groups \cite{dona-kron}
\begin{equation}\label{spin4}
  G = Spin(4) \cong SU(2)_+ \times SU(2)_-.
\end{equation}
The product structure \eq{spin4} also leads to the splitting of the Lie algebra
\begin{equation}\label{lie-spin4}
\mathfrak{g} = spin(4) \cong so(4) \cong su(2)_+ \oplus su(2)_-.
\end{equation}
The splitting of the vector space $\mathfrak{g} = so(4)$ (or $spin(4)$) is defined by the chiral operator
$\gamma_5 = - \gamma_1 \gamma_2 \gamma_3 \gamma_4$ in the Clifford algebra which obeys $\gamma_5^2 =1$.
Therefore one can construct the projection operator
\begin{equation}\label{proj-2}
P_\pm = \frac{1}{2} (1 \pm \gamma_5), \qquad P_\pm^2 = P_\pm, \qquad P_+ + P_- = 1,
\end{equation}
acting on the vector space $\mathfrak{g} = so(4)$.
Using the projection operator \eq{proj-2}, the explicit realization of the splitting \eq{lie-spin4} reads as
\begin{equation}\label{split-so4}
J_{ab} = J_{ab}^+ \oplus J_{ab}^-
\end{equation}
where $J_{ab} = \frac{1}{4} [\gamma_a, \gamma_b ]$ are Lie algebra generators of $so(4)$
and $J_{ab}^\pm  \equiv P_\pm J_{ab}$. Each part in Eq. \eq{split-so4} consists of three $(2 \times 2)$ anti-Hermitian
traceless matrices.
So they can be expanded in the basis of the Pauli matrices $\tau^i \; (i=1,2,3)$ as
\begin{equation}\label{split-su2}
J_{ab}^+ = \frac{i}{2} \eta^i_{ab} \tau^i \in  su(2)_+, \qquad
J_{ab}^- = \frac{i}{2}\bar{ \eta}^i_{ab} \tau^i \in su(2)_-.
\end{equation}
The expansion coefficients in \eq{split-su2}, the so-called 't Hooft symbols \cite{egh-report}, are given by
\begin{equation}\label{t-symbol}
 \eta^i_{ab} = - i \mathrm{Tr} \big( J_{ab}^+ \tau^i \big), \qquad
 \bar{\eta}^i_{ab} = - i \mathrm{Tr} \big( J_{ab}^- \tau^i \big)
\end{equation}
and they satisfy the self-duality relation
\begin{equation}\label{sd-asd-thooft}
\eta^i_{ab} = \frac{1}{2} \varepsilon_{abcd} \eta^i_{cd}, \qquad
\bar{\eta}^i_{ab} = - \frac{1}{2} \varepsilon_{abcd} \bar{\eta}^i_{cd}.
\end{equation}

The decomposition \eq{sd-asd} is that the six-dimensional vector space of two-forms canonically splits
into the sum of three-dimensional vector spaces of self-dual and anti-self-dual two-forms.
Canonical basis elements of self-dual and anti-self-dual two forms are given by
\begin{equation}\label{sd-2form}
  \zeta^i_+ = \frac{1}{2} \eta^i_{ab} e^a \wedge e^b \in  \Omega^2_+, \qquad
  \zeta^i_- = \frac{1}{2} \bar{\eta}^i_{ab}  e^a \wedge e^b \in \Omega^2_-.
\end{equation}
Explicitly, they are given by
\begin{equation}\label{tri-2form}
  \zeta^1_\pm = e^2 \wedge e^3 \pm e^1 \wedge e^4, \quad
  \zeta^2_\pm = e^3 \wedge e^1 \pm e^2 \wedge e^4, \quad
  \zeta^3_\pm = e^1 \wedge e^2 \pm e^3 \wedge e^4.
\end{equation}
They satisfy the Hodge-duality equation
\begin{equation}\label{hodge-bais}
  * \zeta^i_\pm = \pm \zeta^i_\pm.
\end{equation}
and also obey the intersection relation \cite{dona-kron}
\begin{equation}\label{intersection}
\zeta^i_\pm \wedge \zeta^j_\pm = \pm 2 \delta^{ij} d\mu, \qquad
\zeta^i_\pm \wedge \zeta^j_\mp = 0,
\end{equation}
where $d\mu \equiv e^1 \wedge e^2 \wedge e^3 \wedge e^4 = \sqrt{g} d^4 x$ is a volume form on a four-manifold $M$.

Consider a Riemannian manifold $(M,g)$. The metric on $M$ takes the form
\begin{equation}\label{4-metric}
  ds^2 = g_{\mu\nu} (x) dx^\mu dx^\nu = e^a \otimes e^a.
\end{equation}
Using the metric, one can determine the spin connections ${\omega^a}_{b} = {\omega^a}_{b \mu} dx^\mu$
and curvature tensors ${R^a}_{b} = \frac{1}{2} {R^a}_{b \mu \nu} dx^\mu \wedge dx^\nu$
by solving the structure equations \cite{egh-report,nakahara}
\begin{eqnarray} \label{t-free}
  &&  T^a = de^a + {\omega^a}_{b} \wedge e^b =0, \\
  \label{curv-eq}
  && {R^a}_{b} = d {\omega^a}_{b} + {\omega^a}_{c} \wedge {\omega^c}_{b}.
\end{eqnarray}

An important point is that the spin connections ${\omega^a}_{b}$ and Riemann curvature tensors ${R^a}_{b}$ are $so(4)$-valued
one-forms and two-forms in $\Omega^p (M) = \Lambda^p T^* M$, respectively.
Thus one can apply the decompositions in \eq{sd-asd} and \eq{split-so4} to spin connections and curvature tensors \cite{oh-yang,yang-col2,yang-col3}.
The decomposition \eq{lie-spin4} implies that the spin connections can be split into a pair of $SU(2)_+$
and $SU(2)_-$ connections according to the Lie algebra splitting \eq{split-so4}:
\begin{equation}\label{ga-split}
\omega_{ab} = A^{(+)i} \eta^i_{ab}  +  A^{(-)i}\overline{\eta}^i_{ab}.
\end{equation}
The Riemann curvature tensor $R_{ab} = \frac{1}{2} R_{ab\mu\nu} dx^\mu \wedge dx^\nu
= \frac{1}{2} R_{abcd} e^c \wedge e^d \in C^\infty(\mathfrak{g} \otimes \Omega^2)$ is $so(4)$-valued 2-forms.
Thus it is involved with two vector spaces $\mathfrak{g} = so(4)$ in \eq{lie-spin4} and $\Omega^2 (M)$ in \eq{sd-asd}.
Let us apply the decompositions of the vector spaces $\mathfrak{g}$ and $\Omega^2(M)$ to Riemann curvature tensors.
The first decomposition is that the Riemann tensor can be split into a pair of $SU(2)_+$ and $SU(2)_-$ curvatures
according to the Lie algebra splitting \eq{split-so4}:
\begin{equation}\label{split-riemann}
  R_{ab} = F^{(+)i} \eta^i_{ab}  + F^{(-)i} \bar{\eta}^i_{ab},
\end{equation}
where $SU(2)_\pm$ curvatures are 2-forms on $M$ defined by
\begin{eqnarray}\label{sd-2f}
  F^{(\pm)i} &=& \frac{1}{2} F^{(\pm)i}_{cd} e^c \wedge e^d \nonumber \\
             &=& dA^{(\pm)i} - \varepsilon^{ijk} A^{(\pm)j} \wedge A^{(\pm)k}.
\end{eqnarray}
The second decomposition is that the $SU(2)_\pm$ curvatures in Eq. \eq{sd-2f} can be decomposed as
\begin{equation}\label{decop-2f}
 F^{(+)i} =  f^{ij}_{(++)}\zeta^j_{+} + f^{ij}_{(+-)}\zeta^j_{-}, \qquad
  F^{(-)i} =  f^{ij}_{(-+)}\zeta^j_{+} + f^{ij}_{(--)}\zeta^j_{-}.
\end{equation}
Using the intersection pairing \eq{intersection}, the coefficients are determined as
\begin{eqnarray} \label{coeff-1}
    && f^{ij}_{(\pm \pm)} d \mu = \pm \frac{1}{4} \left( F^{(\pm)i} \wedge \zeta^j_{\pm}
    + F^{(\pm)j} \wedge \zeta^i_{\pm} \right), \\
    \label{coeff-2}
    && f^{ij}_{(\pm \mp)} d \mu = \mp \frac{1}{4} \left( F^{(\pm)i} \wedge \zeta^j_{\mp}
    - F^{(\mp)j} \wedge \zeta^i_{\pm} \right).
\end{eqnarray}
Combining the two decompositions \eq{split-riemann} and \eq{decop-2f} leads to an irreducible decomposition of
the general Riemann curvature tensor \cite{oh-yang,yang-col2}:
\begin{equation}\label{decom-riem}
  R_{abcd} = f^{ij}_{(++)}\eta^i_{ab}\eta^j_{cd} + f^{ij}_{(+-)}\eta^i_{ab} \bar{\eta}^j_{cd}
  +  f^{ij}_{(-+)} \bar{\eta}^i_{ab} \eta^j_{cd} + f^{ij}_{(--)} \bar{\eta}^i_{ab}\bar{\eta}^j_{cd}.
\end{equation}
Recently, the decomposition \eq{decom-riem} has been used significantly in \cite{chy-grg}
to study the scalar invariants of the curvature tensor.

Now we impose the first Bianchi identity
\begin{equation}\label{1-bianchi}
R_{abcd} + R_{acdb} + R_{adbc} = 0
\end{equation}
to the decomposition \eq{decom-riem}.
The Bianchi identity, being totally 16 conditions, can equivalently be stated as
\begin{eqnarray} \label{bianchi-symm}
&& R_{abcd} = R_{cdab}, \\
\label{bianchi-pseudo}
&& \varepsilon^{abcd} R_{abcd} = 0.
\end{eqnarray}
Eq. (\ref{bianchi-symm}) requires the symmetry property
\begin{equation}\label{symm-fij}
f^{ij}_{(++)} = f^{ji}_{(++)}, \quad f^{ij}_{(--)}=f^{ji}_{(--)},
\quad  f^{ij}_{(+-)} = f^{ji}_{(-+)}.
\end{equation}
whereas Eq. (\ref{bianchi-pseudo}) imposes an additional constraint
\begin{equation}\label{1-more}
  f^{ij}_{(++)} \delta^{ij} = f^{ij}_{(--)} \delta^{ij}.
\end{equation}

The ``trace-free part" of the Riemann curvature tensor is called the Weyl tensor defined by \cite{egh-report,nakahara}
\begin{equation}\label{weyl}
  W_{abcd} = R_{abcd} - \frac{1}{2} (\delta_{ac} R_{bd} - \delta_{ad} R_{bc} - \delta_{bc} R_{ad} + \delta_{bd} R_{ac} )
  + \frac{1}{6} (\delta_{ac} \delta_{bd} - \delta_{ad} \delta_{bc}) R.
\end{equation}
The Weyl tensor shares all the symmetry structures of the curvature tensor and all its traces with the
metric are zero. Then the Weyl tensor can be decomposed as \cite{oh-yang,yang-col2}
\begin{eqnarray}\label{weyl-riem}
  W_{abcd}
  = {\widetilde f}^{ij}_{(++)} \eta^i_{ab}\eta^j_{cd} + {\widetilde f}^{ij}_{(--)} \bar{\eta}^i_{ab}\bar{\eta}^j_{cd},
\end{eqnarray}
where ${\widetilde f}^{ij}_{(\pm\pm)} \equiv f^{ij}_{(\pm\pm)} - \frac{1}{3} \delta^{ij} \big(  f^{kl}_{(\pm\pm)} \delta^{kl} \big)$
are symmetric, traceless $3 \times 3$ matrices.

The decomposition \eq{decom-riem} can be applied to the Ricci tensor $R_{ab} \equiv R_{acbc}$ and
the Ricci scalar $R \equiv R_{aa}$ to yield
\begin{eqnarray}\label{ricci}
  R_{ab} &=& \big(f^{ij}_{(++)} \delta^{ij} + f^{ij}_{(--)} \delta^{ij} \big) \delta_{ab}
  + 2 f^{ij}_{(+-)}\eta^i_{ac} \bar{\eta}^j_{bc}, \nonumber \\
  R &=& 4 \big(f^{ij}_{(++)} \delta^{ij} + f^{ij}_{(--)} \delta^{ij} \big).
\end{eqnarray}
For Einstein manifolds satisfying the equations, $R_{ab} = \lambda \delta_{ab}$, with $\lambda$ a cosmological constant, one can show \cite{oh-yang,yang-col2} that
\begin{equation} \label{einstein}
  f^{ij}_{(+-)} = f^{ji}_{(-+)} = 0.
\end{equation}
As a result, the curvature tensor for Einstein manifolds takes the simple form
\begin{equation}\label{em-riem}
  R_{abcd} = f^{ij}_{(++)}\eta^i_{ab}\eta^j_{cd} + f^{ij}_{(--)} \bar{\eta}^i_{ab}\bar{\eta}^j_{cd}.
\end{equation}
Therefore, we arrive at the following Lemma, proved in [3]:

$\mathbf{Lemma \; 2.1}$: The Einstein equations, $R_{ab} = \lambda \delta_{ab}, \; \lambda \in \mathbb{R}$,
for an orientable Riemannian manifold $(M, g)$ are equivalent to the self-duality equations
of Yang-Mills instantons,
\begin{equation}\label{sde}
F^{(\pm)i} = \pm * F^{(\pm)i},
\end{equation}
where $F^{(\pm)i} \wedge \zeta^i_\pm = \pm \lambda d \mu$.

\noindent It should be pointed out that the equivalence of the equations does not mean the equivalence of solution spaces because
the Einstein equations are the second-order differential equations of an Einstein metric while the self-duality equations are
the first-order differential equations of an $SU(2)$ connection. For example, if an Einstein metric is given,
then $SU(2)$ connections are determined up to gauge transformations by solving Eq. \eq{ga-split} \cite{yang-col1}.
But the converse is not true.
In order to determine the metric from given $SU(2)$ connections which determine the spin connection via Eq. \eq{ga-split},
it is necessary to solve the torsion-free equation \eq{t-free}. Moreover, it is necessary to specify the pair of instantons and anti-instantons to determine the spin connections. Therefore the moduli space of Yang-Mills instantons
is not directly related to the moduli space of Einstein metrics. Nevertheless, since non-Einstein metrics
imply non-instanton connections, the Lemma $\mathbf{2.1}$ will be very useful for studying the moduli space
of Einstein structures.

The above understanding of orientable four-manifolds comes via the geometry of connections or gauge fields
where Riemannian manifolds are described as a gauge theory of $SU(2)_\pm$ connections in \eq{ga-split}
and their curvatures in \eq{split-riemann}.
Therefore Riemannian manifolds encode a topological information in the form of Yang-Mills instantons of $SU(2)_\pm$ gauge fields.
Hence it is natural to expect that the topological invariants of a Riemannian manifold $(M,g)$ will be determined
by the configuration of $SU(2)_\pm$ Yang-Mills instantons.
For a general closed Riemannian manifold $M$ without boundaries, the Euler characteristic $\chi(M)$ and
the Hirzebruch signature $\tau(M)$ are defined by \cite{besse,gh-cmp1979,egh-report}
\begin{eqnarray} \label{top-euler}
\chi(M) &=& \frac{1}{32 \pi^2} \int_M \varepsilon^{abcd} R_{ab} \wedge R_{cd}, \\
\label{top-hirze}
\tau(M) &=& \frac{1}{24 \pi^2} \int_M  R_{ab} \wedge R_{ab}.
\end{eqnarray}
The topological invariants can be expressed in terms of $SU(2)_\pm$ connections
using the decompositions \eq{split-riemann} and \eq{decop-2f} \cite{oh-yang,yang-col3,yang-col4}
\begin{eqnarray} \label{dec-euler}
\chi(M) &=& \frac{1}{4 \pi^2} \int_M \big( F^{(+)i} \wedge F^{(+)i} - F^{(-)i} \wedge F^{(-)i} \big), \xx
&=& \frac{1}{2 \pi^2} \int_M \left( (f^{ij}_{(++)})^2 + (f^{ij}_{(--)})^2 - 2 (f^{ij}_{(+-)})^2 \right) d\mu, \\
\label{dec-hirze}
\tau(M) &=& \frac{1}{6 \pi^2} \int_M \big( F^{(+)i} \wedge F^{(+)i} + F^{(+)i} \wedge F^{(+)i} \big), \xx
&=& \frac{1}{3 \pi^2} \int_M \left( (f^{ij}_{(++)})^2 - (f^{ij}_{(--)})^2 \right) d\mu,
\end{eqnarray}
where we used the intersection relation \eq{intersection}.

\begin{figure}
\centering
\includegraphics[width=0.6\textwidth]{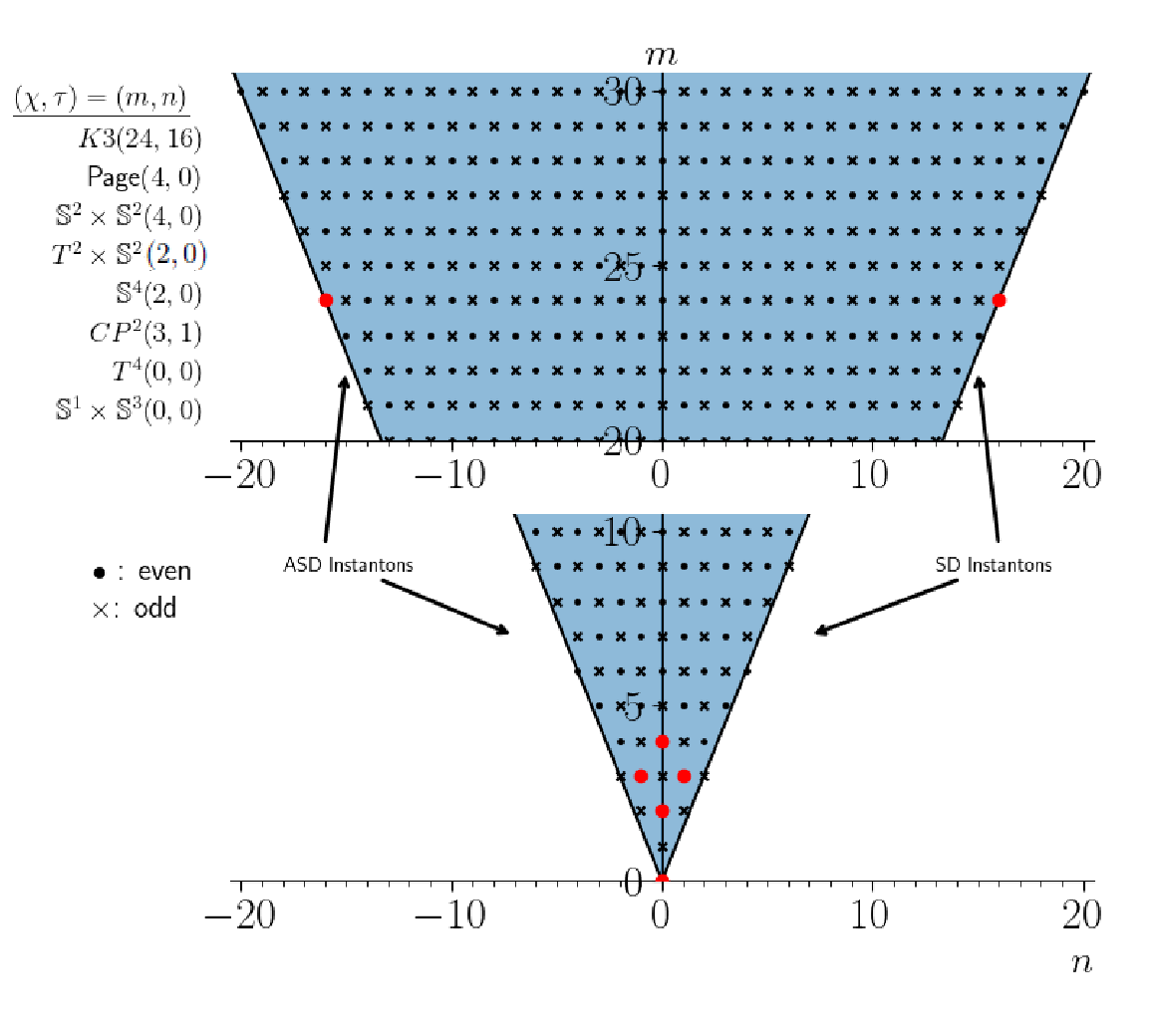}
\caption{\label{top-num} Topological numbers of closed Einstein manifolds \cite{yang-col4}}
\end{figure}

An Einstein manifold has curvature tensors given by \eq{em-riem}.
In this case, the Euler characteristic $\chi(M)$ is determined by the sum of instantons and anti-instantons
whereas the Hirzebruch signature $\tau(M)$ is their difference. Then it becomes easy to verify
the important inequalities for the topological invariants. The first inequality is $\chi(M)\geq 0$ with equality
only if $f^{ij}_{(++)} = f^{ij}_{(--)} = 0$, i.e., $M$ is flat ($\mathbf{6.32}$ in \cite{besse}
and Sect. 10.4 in \cite{egh-report}).
The second inequality is the Hitchin-Thorpe inequality \cite{hitchin1974,pps-blms} stating that
\begin{equation}\label{hitchin-thorpe}
  \chi (M) \pm \frac{3}{2} \tau (M) = \frac{1}{\pi^2} \int_M  \big(f^{ij}_{(\pm\pm)} \big)^2 d\mu \geq 0
\end{equation}
where the equality holds only if $f^{ij}_{(++)} = 0$ or $f^{ij}_{(--)} = 0$, i.e., $M$ is half-flat
(a gravitational instanton). A direct consequence of the first inequality is that $\mathbb{S}^1 \times \mathbb{S}^3$
does not admit Einstein metrics \cite{besse}. On the other hand, the second inequality leads to the result
that $n \mathbb{C}\mathbb{P}^2 \; (n \geq 4)$ is a simply connected compact manifold
which does not carry an Einstein metric \cite{hitchin1974}, since $\chi=n+2$ and $\tau = n$.
For the rational surfaces $M_k = \mathbb{C}\mathbb{P}^2 \sharp k \overline{\mathbb{C}\mathbb{P}}^2$,
one has $\chi(M_k) = 3 + k$ and $\tau(M_k) = 1 - k$, so the Hitchin-Thorpe inequality implies that $M_k$
does not admit an Einstein metric if $k \geq 9$ \cite{rev-anderson}.
The four-manifold $M_k$ does admit an Einstein metric only for $k \leq 8$.
Some examples of four-dimensional compact Einstein
manifolds are shown up in Fig. \ref{top-num}.

If an Einstein manifold with positive scalar curvature is further assumed to have a nonnegative sectional curvature,
more refined inequalities can be stated \cite{hitchin1974}. A compact Einstein manifold $M$ with
nonnegative sectional curvature must have its Euler characteristic $\chi(M)$ bounded by $1 \leq \chi(M) \leq 9$ and
its Hirzebruch signature $\tau(M)$ should satisfy the inequality
\begin{equation}\label{stineq-hitchin}
  |\tau (M)| \leq \left(\frac{2}{3} \right)^{\frac{3}{2}} \chi(M).
\end{equation}
We will prove this inequality in appendix A.
Since $\left(\frac{2}{3} \right)^{\frac{3}{2}} \approx 0.544$ is an irrational number,
the equality can only occur if $M$ is flat. Then we see that $|\tau (M)| \leq \left(\frac{2}{3} \right)^{\frac{3}{2}}
\times 9 \approx 4.899$, i.e., $|\tau (M)| \leq 4$.
These results show that most four-manifolds cannot carry any Einstein structure with positive
or nonnegative sectional curvature. This is the reason why the red dots in Fig. \ref{top-num}
are mostly located in small topological numbers except the $K3$ surface which is a half-flat manifold.

There is an important theorem on compact complex surfaces $(M, J)$ \cite{tian-yau,tian-im} where $J$ is the complex structure on $M$:

$\mathbf{Theorem \; 2.2}$: A compact complex surface $(M, J)$ admits a K\"ahler-Einstein metric with $\lambda > 0$
if and only if the first Chern class $c_1(M)$ is positive and the Lie algebra of holomorphic vector fields is reductive.
This occurs exactly on $\mathbb{C}\mathbb{P}^2, \; \mathbb{S}^2 \times \mathbb{S}^2$,
or $\mathbb{C}\mathbb{P}^2 \sharp k \overline{\mathbb{C}\mathbb{P}}^2 \; (3 \leq k \leq 8)$ which is the blow-up of $\mathbb{C}\mathbb{P}^2$ at general $k$ points.

\noindent Thus the Page space $\mathbb{C}\mathbb{P}^2 \sharp \, \overline{\mathbb{C}\mathbb{P}}^2$, which will be considered below,
is not a K\"ahler manifold. Another non-K\"ahler metric on $\mathbb{C}\mathbb{P}^2 \sharp 2 \overline{\mathbb{C}\mathbb{P}}^2$
was constructed in \cite{lebrun-jams}. Note that there are relatively few K\"ahler-Einstein metrics when $\lambda \geq 0$;
most all have $\lambda < 0 $, in analogy to the case of Riemann surfaces.

It will be instructive to analyze the Einstein structure of the Page space \cite{page} in detail
since it is the first example of nonhomogeneous Einstein manifolds and its explicit metric is known.
The Page space $\mathbb{C}\mathbb{P}^2 \sharp \, \overline{\mathbb{C}\mathbb{P}}^2$ can be constructed
by the blow-up of $\mathbb{C}\mathbb{P}^2$ at a point $p \in \mathbb{C}\mathbb{P}^2$  \cite{besse}.
The metric on the Page space takes on a complicated form
\begin{eqnarray}\label{page}
  ds^2 &=& e^a \otimes e^a \\
  &=& \frac{3}{\lambda} \frac{1+\nu^2}{3+\nu^2} \left( U^{-1} d\chi^2
  + \frac{U}{4} \sin^2 \chi (d\psi+ \cos \theta d\phi)^2
  + \frac{1-\nu^2 \cos^2 \chi}{4\nu} (d\theta^2 + \sin^2 \theta d\phi^2) \right), \nonumber
\end{eqnarray}
where
$$ U(\chi) = 1 - \frac{2 \nu^2 \sin^2 \chi}{(3+ \nu^2)(1-\nu^2 \cos^2 \chi)} $$
and $\nu$ is the positive root of $\nu^4 + 4 \nu^3 - 6\nu^2 + 12\nu - 3=0$ which works out to be $0.2817$.
The spin connections for the Page metric \eq{page} can be determined by solving the torsion-free condition \eq{t-free} as
\begin{eqnarray}
&& \omega_{12} = - \sqrt{\frac{4\nu}{1-\nu^2 \cos^2 \chi}} \cot \theta e^2 + f e^3,
\quad \omega_{13} = f e^2, \quad \omega_{23} = - f e^1,   \\
&& \omega_{14} = g e^1,  \quad  \omega_{24} = g e^2, \quad
\omega_{34} = \left( \sqrt{U} \cot \chi + \frac{2 (\nu^2 -1)}{U (3+\nu^2) (1-\nu^2 \cos^2 \chi)} g  \right) e^3,\nonumber
\end{eqnarray}
where
$$ f= \frac{\sqrt{U} \nu}{1-\nu^2 \cos^2 \chi} \sin \chi,  \qquad g =  f \nu \cos \chi.$$
Here we chose $\lambda = 3 \frac{1+\nu^2}{3+\nu^2}$ for simplicity.
The corresponding $SU(2)$ connections are given by \cite{hkky}
\begin{eqnarray*} \label{page-gauge}
  A^{(\pm)1} &=& \frac{1}{2} (-f \pm g) e^1, \qquad A^{(\pm)2} = \frac{1}{2} (-f \pm g) e^2, \\
  A^{(\pm)3} &=& \frac{1}{2} \left( - \sqrt{\frac{4\nu}{1-\nu^2 \cos^2 \chi}} \cot \theta e^2 + f e^3
  \pm \Big( \sqrt{U} \cot \chi + \frac{2 (\nu^2 -1)}{U (3+\nu^2) (1-\nu^2 \cos^2 \chi)} g \Big) e^3 \right).
\end{eqnarray*}
The corresponding $SU(2)$ curvatures have the form, $f^{ij}_{(\pm\pm)} = \pm \psi^i_{(\pm)} \delta^{ij}$, i.e.,
\begin{eqnarray}\label{page-f}
  F^{(\pm)1} = \pm \psi^1_{(\pm)} \zeta^1_\pm, \qquad F^{(\pm)2} = \pm \psi^2_{(\pm)} \zeta^2_\pm, \qquad
  F^{(\pm)1} = \pm \psi^3_{(\pm)} \zeta^3_\pm,
\end{eqnarray}
where $ \psi^1_{(\pm)} =  \psi^2_{(\pm)}$ and
\begin{eqnarray} \label{page-comp}
  \psi^1_{(\pm)} &=& \frac{U}{2} \frac{\nu(1-\nu^2) \cos \chi \mp \nu^2 (\cos^2 \chi - \sin^2 \chi)
  \pm \nu^4 \cos^4 \chi}{(1-\nu^2 \cos^2 \chi)^2} \xx
  && - \frac{\nu^3 (1-\nu^2) \cos \chi \sin^2 \chi}{(3+\nu^2)(1 - \nu^2 \cos^2 \chi)^3} (1 \mp \nu\cos \chi), \xx
  \psi^3_{(\pm)} &=& - \frac{\nu}{2(1 - \nu^2 \cos^2 \chi)}
  \left( U' \sin \chi + \frac{2(1-\nu^2) U \cos \chi}{1 - \nu^2 \cos^2 \chi} \right) \xx
&& \pm \frac{\nu}{2(1 - \nu^2 \cos^2 \chi)}
  \left( 4-  \frac{\nu U \sin^2 \chi (3+ \nu^2 \cos^2 \chi)}{1 - \nu^2 \cos^2 \chi} \right).
\end{eqnarray}
The result \eq{page-f} shows that the $SU(2)_+ \; (SU(2)_-)$ curvatures of the Page metric \eq{page}
are self-dual (anti-self-dual), so the Page space must be an Einstein manifold according to the Lemma $\mathbf{2.1}$.
One can explicitly check that the Ricci tensor is given by (after recovering the factor $\frac{3}{\lambda} \frac{1+\nu^2}{3+\nu^2}$)
\begin{align}
R_{ab} = \lambda \, \text{diag} \left( c, c, 1, 1 \right),
\end{align}
where $c = \frac{\nu^4-8 \nu^3 + 9 \nu^2 - 24 \nu +3 \left(\nu^2+1\right) \nu^2 \cos 2 \chi }{3 \nu^4-3 \nu^2 - 6 + 3
\left(\nu ^2 + 1\right) \nu^2 \cos 2 \chi }$. Using the quartic equation $\nu^4 + 4 \nu^3 - 6\nu^2 + 12\nu - 3=0$,
one can show that $c=1$. We show in appendix A that (anti-)self-dual Weyl tensors have only two distinct eigenvalues at each point.
So, according to Proposition 5 in \cite{v21}, the Page metric \eq{page} is locally conformal to a K\"ahler metric
although it is not a K\"ahler metric in itself.

The explicit calculation using the result \eq{page-f} gives us $\chi(M) = 4$ and $\tau(M)=0$
for the Page space \cite{gh-cmp1979} as shown up in Fig. \ref{top-num}.
It should be $\tau(M)=0$ if $M$ is the manifold obtained by blowing up one point
in $\mathbb{C}\mathbb{P}^2$ (see $\mathbf{11.82 \; Example}$ in \cite{besse}).

\section{Squashed Four-Spheres}

Let us define the instanton density as $\rho^{(\pm)}_n = - \frac{1}{4} \mathrm{Tr} F^{(\pm)}_{\alpha\beta}F^{(\pm) \alpha\beta}$,
which describes a curvature distribution over $M$.
The following formula is useful;
\begin{equation}\label{density}
  F^{(\pm)i} \wedge F^{(\pm)i} = \pm 2 \rho^{(\pm)}_n d \mu.
\end{equation}
If $M$ is a homogeneous Einstein manifold, the left-hand side in Eq. \eq{density} must be a constant multiple of the volume form.
In other words, the instanton density \eq{density} characterizes the homogeneity of an underlying Einstein manifold.
One can check that the Page space \eq{page} has a nontrivial distribution of the instanton density \cite{hkky},
so confirm the nonhomogeneity. Using the definition \eq{density},
the topological invariants in \eq{dec-euler} and \eq{dec-hirze} can be written as
\begin{eqnarray} \label{den-euler}
\chi(M) &=& \int_M  \rho_\chi (M) d \mu, \\
\label{den-hirze}
\tau(M) &=& \int_M \rho_\tau (M) d \mu,
\end{eqnarray}
where
\begin{eqnarray}
\label{def-de}
\rho_\chi (M) &=& \frac{1}{2 \pi^2} ( \rho^{(+)}_n  + \rho^{(-)}_n ),  \\
\label{def-dh}
\rho_\tau (M) &=& \frac{1}{3 \pi^2} ( \rho^{(+)}_n  - \rho^{(-)}_n ).
\end{eqnarray}
Thus $\rho_\chi (M) \pm \frac{3}{2} \rho_\tau (M) = \frac{1}{\pi^2} \rho^{(\pm)}_n \geq 0$.

We now examine another nonhomogeneous space obtained by deforming a round sphere $\mathbb{S}^4$
whose deformation is parameterized by three axis-length parameters $(r, l, \tilde{l})$.
The four-dimensional ellipsoids are defined by
\begin{equation}\label{4d-ell}
  \frac{x_0^2}{r^2} + \frac{x_1^2 + x_2^2}{l^2} + \frac{x_3^2 + x_4^2}{\tilde{l}^2} =1.
\end{equation}
The metric of squashed 4-sphere is given by $ds^2 = e^a \otimes e^a$ with the vierbeins \cite{squashed-s4}
\begin{equation} \label{ss4-vierbein}
  \left\{ e^1 = \sin \rho \tilde{\sigma}^1, \; e^2 = \sin \rho \tilde{\sigma}^2, \;
e^3 = \sin \rho \tilde{\sigma}^3 + h d\rho, \; e^4 = g d\rho \right\},
\end{equation}
where
\begin{eqnarray*}
  f &\equiv & \sqrt{l^2 \sin^2 \theta + \tilde{l}^2 \cos^2 \theta}, \\
  g &\equiv & \sqrt{r^2 \sin^2 \rho + l^2 \tilde{l}^2 f^{-2} \cos^2 \rho}, \\
  h &\equiv & \frac{\tilde{l}^2 -l^2}{f} \cos \rho \sin \theta \cos \theta,
\end{eqnarray*}
and $\tilde{\sigma}^i \, (i=1,2,3)$ are dreibeins of the three-dimensional ellipsoid
given in polar coordinates $(\varphi, \chi, \theta)$ by
\begin{equation*}
  \tilde{\sigma}^1 = l \cos \theta d \varphi, \quad
  \tilde{\sigma}^2 = \tilde{l} \sin \theta d \chi, \quad \tilde{\sigma}^3 = f d \theta.
\end{equation*}
The ranges of coordinates are $\rho \in [0, \pi], \; \theta \in [0, \frac{\pi}{2}]$ and
$\varphi, \chi \in [0, 2\pi]$.

It is straightforward to calculate the $SU(2)$ connections
\begin{eqnarray} \label{hhaho-gauge}
  A^{(\pm)1} &=& \pm \frac{\tilde{l}^2}{2 gf^2} \cot \rho e^1 + \frac{\cot \theta}{2 f \sin \rho}e^2, \xx
  A^{(\pm)2} &=& \pm \frac{l^2}{2 gf^2} \cot \rho e^2 + \frac{\tan \theta}{2 f \sin \rho}e^1, \\
  A^{(\pm)3} &=& \pm \frac{l^2 \tilde{l}^2}{2 gf^4} \cot \rho \left( e^3 - \frac{h}{g} e^4 \right), \nonumber
\end{eqnarray}
and their curvatures
\begin{eqnarray}\label{ss4-f}
  F^{(\pm)1} &=& \pm \frac{\tilde{l}^2}{2 f^3 \sin^2 \rho} \left(
  \frac{h}{g^3 f^3} (2 g^2 f^2 - l^2 \tilde{l}^2 \cos^2 \rho )
 - \frac{(\tilde{l}^2 - l^2)}{2} \psi \sin 2 \theta \right) e^3 \wedge e^1 \xx
 && \pm \frac{\tilde{l}^2}{2 g f^2 \sin^2 \rho} \left(
  \frac{r^2}{g^3} \sin^2 \rho + \frac{h^2}{g^3 f^4} (2 g^2 f^2 - l^2 \tilde{l}^2 \cos^2 \rho )
 - \frac{(\tilde{l}^2 - l^2) h}{2f} \psi \sin 2 \theta \right) e^1 \wedge e^4 \xx
 && + \frac{l^2 r^2}{2 f^4 g^2} e^2 \wedge e^3 - \frac{l^2 r^2 h}{2 f^4 g^3} e^2 \wedge e^4  \xx
  F^{(\pm)2} &=& \pm \frac{l^2}{2 f^3 \sin^2 \rho} \left(
  - \frac{h}{g^3 f^3} (2 g^2 f^2 - l^2 \tilde{l}^2 \cos^2 \rho )
 + \frac{(\tilde{l}^2 - l^2)}{2} \psi \sin 2 \theta \right) e^2 \wedge e^3 \\
 && \pm \frac{l^2}{2 g f^2 \sin^2 \rho} \left(
  \frac{r^2}{g^3} \sin^2 \rho + \frac{h^2}{g^3 f^4} (2 g^2 f^2 - l^2 \tilde{l}^2 \cos^2 \rho )
 - \frac{(\tilde{l}^2 - l^2) h}{2f} \psi \sin 2 \theta \right) e^2 \wedge e^4 \xx
 && + \frac{\tilde{l}^2 r^2}{2 f^4 g^2} e^3 \wedge e^1 + \frac{\tilde{l}^2 r^2 h}{2 f^4 g^3} e^1 \wedge e^4  \xx
  F^{(\pm)3} &=& \frac{r^2}{2f^2 g^2} \left( e^1 \wedge e^2 \pm \frac{l^2 \tilde{l}^2}{f^2 g^2} e^3 \wedge e^4 \right), \nonumber
\end{eqnarray}
where $\psi \equiv \frac{\cos \rho}{ g f^2}$.

The squashed four-sphere is not an Einstein manifold. According to the Lemma $\mathbf{2.1}$,
$SU(2)$ connections constructed from a four-manifold $M$ become self-dual or anti-self-dual only
if $M$ is Einstein, i.e., obeying the vacuum Einstein equations, $R_{\mu\nu} = \lambda g_{\mu\nu}$.
One can check that the $SU(2)$ curvatures in Eq. \eq{ss4-f} do not satisfy the self-duality \eq{sde}.
The instanton density for the metric \eq{ss4-vierbein} exhibits a nontrivial distribution
\begin{eqnarray} \label{inst-dens}
\rho^{(\pm)}_n = \frac{3}{4} \frac{l^2 \tilde{l}^2 r^4}{f^6 g^6}.
\end{eqnarray}
Therefore the squashed $\mathbb{S}^4$ is nonhomogeneous as expected.
When $l = \tilde{l} = r$, it becomes the round four-sphere $\mathbb{S}^4$ which is an Einstein manifold.
In this case, the $SU(2)$ connections in \eq{sqs4-gauge} become instanton connections\footnote{The connections
on $\mathbb{S}^4$ can be mapped to the gauge fields of a self-dual SU(2) instanton on $\mathbb{R}^4$
by using the $\mathbb{R}^4$-coordinates for the left-invariant one-forms $\sigma^i$ on $\mathbb{S}^3$ \cite{hkky}.} and
the instanton density \eq{inst-dens} becomes constant, i.e., $F^{(\pm)i} \wedge F^{(\pm)i} = \pm \frac{3}{2 r^4} d \mu$.

The deformation of Einstein structures is captureed by the coefficients $f^{ij}_{(\pm \mp)}$ which can be computed
by Eq. \eq{coeff-2}:
\begin{eqnarray} \label{def-coeff}
f^{11}_{(\pm \mp)} &=& - \frac{\tilde{l}^2}{4 g f^2 \sin^2 \rho} \left(
  \frac{r^2}{g^3} \sin^2 \rho + \frac{h^2}{g^3 f^4} (2 g^2 f^2 - l^2 \tilde{l}^2 \cos^2 \rho )
 - \frac{(\tilde{l}^2 - l^2) h}{2f} \psi \sin 2 \theta \right) + \frac{l^2 r^2}{4 f^4 g^2},  \xx
f^{12}_{(\pm \mp)} &=& \pm \frac{(l^2 + \tilde{l}^2)}{8 f^3 \sin^2 \rho} \left(
  \frac{h}{g^3 f^3} (2 g^2 f^2 - l^2 \tilde{l}^2 \cos^2 \rho )
 - \frac{(\tilde{l}^2 - l^2)}{2} \psi \sin 2 \theta + \frac{r^2 h}{f g^3} \sin^2 \rho \right), \xx
f^{22}_{(\pm \mp)} &=& - \frac{l^2}{4 g f^2 \sin^2 \rho} \left(
  \frac{r^2}{g^3} \sin^2 \rho + \frac{h^2}{g^3 f^4} (2 g^2 f^2 - l^2 \tilde{l}^2 \cos^2 \rho )
 - \frac{(\tilde{l}^2 - l^2) h}{2f} \psi \sin 2 \theta \right) + \frac{\tilde{l}^2 r^2}{4 f^4 g^2},  \xx
f^{33}_{(\pm \mp)} &=& \frac{r^2}{4f^2 g^2} \left( 1- \frac{l^2 \tilde{l}^2}{f^2 g^2} \right),
\end{eqnarray}
and $f^{13}_{(\pm \mp)} = f^{23}_{(\pm \mp)} = f^{31}_{(\pm \mp)} = f^{32}_{(\pm \mp)} = 0$.
If the coefficients in Eq. \eq{def-coeff} vanish, the metric \eq{ss4-vierbein} now becomes an Einstein metric.
Therefore a question is whether there is another Einstein manifold obtained
by continuous deformations of the round four-sphere.
Note that the deformations in \eq{4d-ell} correspond to the deformation of Einstein structures within a fixed conformal class.
Hence, the rigidity of the Einstein structure \cite{besse,daga-yang} implies that
there is no solution of the vanishing condition except as $l = \tilde{l} = r$.
It is easy to see that the coefficients in Eq. \eq{def-coeff} vanish only if $r = l = \tilde{l}$.

The above result implies that there is no continuous deformation of Einstein structures.
Therefore we confirm that the moduli space of Einstein structures is an isolated point.
Now we will generalize the result by considering the family of squashed 4-spheres.
The metric is defined by  \cite{squashed-gs4}
\begin{equation}\label{sq-s4}
  ds^2 = dr^2 + \frac{f(r)^2}{4} (d\theta^2 + \sin^2 \theta d \phi^2)
  + \frac{g(r)^2}{4} (d\psi + \cos \theta d \phi)^2
\end{equation}
and vierbein one-forms are
\begin{equation*}\label{sps4-vierbein}
  e^a = \left\{ \frac{f(r)}{2} d\theta, \frac{f(r)}{2} \sin \theta d \phi,
  \frac{g(r)}{2} (d\psi + \cos \theta d \phi), dr \right\}
\end{equation*}
where $f(r)$ and $g(r)$ are smooth functions of $r$.\footnote{For the metric \eq{sq-s4},
$0 \leq r \leq  \pi, \;0 \leq \theta \leq \pi, \; 0 \leq \phi \leq 2 \pi, \; 0 \leq \psi \leq 4 \pi$.
See Eq. (55) in \cite{squashed-gs4}.}
The above metric has $SU(2) \times U(1)$ isometry \cite{mpps-duke}.
The round four-sphere $\mathbb{S}^4$ corresponds to $f(r) = g(r) = \sin r$.
Note that deformations described by the functions $f(r)$ and $g(r)$ can also include
a deformation of conformal structures of an Einstein metric \cite{pps-im1996,taka-tjm,bonn-cqg}.
However, we will show that there exists only a trivial deformation, at most only changing the size of the sphere,
if the Einstein structure is maintained.

It is straightforward to calculate the $SU(2)$ connections from the metric \eq{sq-s4}
\begin{eqnarray} \label{sqs4-gauge}
  A^{(\pm)1} &=& \frac{1}{2} \left(- \frac{g}{f^2} \pm \frac{f'}{f} \right) e^1,
  \qquad A^{(\pm)2} = \frac{1}{2} \left(- \frac{g}{f^2} \pm \frac{f'}{f} \right) e^2, \xx
  A^{(\pm)3} &=& - \frac{\cot \theta}{f} e^2 + \frac{1}{2} \left( \frac{g}{f^2} \pm \frac{g'}{g} \right) e^3
\end{eqnarray}
and their curvatures
\begin{eqnarray} \label{sqs4-f}
  F^{(\pm)1} &=& \frac{1}{2f^3} (f g'-f'g) (e^1 \wedge e^4 \pm e^2 \wedge e^3)
  \mp \frac{f''}{2f} e^1 \wedge e^4 - \frac{1}{2f} \left( \frac{f'g'}{g} - \frac{g^2}{f^3}\right) e^2 \wedge e^3, \xx
  F^{(\pm)2} &=& \frac{1}{2f^3} (f g'-f'g) (e^2 \wedge e^4 \pm e^3 \wedge e^1)
  \mp \frac{f''}{2f} e^2 \wedge e^4 - \frac{1}{2f} \left( \frac{f'g'}{g} - \frac{g^2}{f^3}\right) e^3 \wedge e^1, \\
  F^{(\pm)3} &=& - \frac{1}{f^3} (f g'-f'g)  (e^3 \wedge e^4 \pm e^1 \wedge e^2 ) \mp \frac{g''}{2g} e^3 \wedge e^4
  - \frac{1}{2f^4} (f^2 f'^2 + 3g^2 - 4f^2) e^1 \wedge e^2. \nonumber
\end{eqnarray}
The instanton density on the four-manifold $(M, g)$ is given by
\begin{eqnarray}
F^{(\pm)i} \wedge F^{(\pm)i} &=& \left[ -\frac{2}{f^4} (fg'-f' g) f''
+ \frac{2 g''}{f^3 g} (fg'-f' g) \right. \xx
&& \;\; \left. \pm \left( \frac{3}{f^6} (fg'-f'g)^2
+ \Big(\frac{f''}{f} \Big)^2 + \frac{1}{2} \Big( \frac{g''}{g} \Big)^2 \right) \right] d \mu \\
&=& \left[ \frac{1}{f^4 g} \frac{d}{dr}(fg'-f' g)^2 \pm \left( \frac{3}{f^6} (fg'-f'g)^2
+ \Big(\frac{f''}{f} \Big)^2 + \frac{1}{2} \Big( \frac{g''}{g} \Big)^2 \right) \right] d \mu. \nonumber
\end{eqnarray}
For the round four-sphere $\mathbb{S}^4$, $f(r) = g(r) = \sin r$ and $F^{(\pm)i} \wedge F^{(\pm)i} = \pm \frac{3}{2} d \mu$.
But the squashed $\mathbb{S}^4$ has $f(r) \neq g(r)$ and the instanton density is not uniform in general
because the underlying manifold is not homogeneous.

One can see that the self-duality \eq{sde} is satisfied if the functions $f(r)$ and $g(r)$ obey
the differential equations
\begin{equation}\label{sqs4-deq}
  \frac{f''}{f} = \frac{f'g'}{fg} - \frac{g^2}{f^4}, \qquad
  \frac{g''}{g} = \frac{f'^2 - 4}{f^2} + \frac{3g^2}{f^4}.
\end{equation}
Then the squashed $\mathbb{S}^4$ becomes an Einstein manifold according to the Lemma $\mathbf{2.1}$.
Consider infinitesimal deformations from the round sphere:
\begin{equation}\label{inf-def}
f(r) = \sin r + \epsilon p(r), \qquad g(r) = \sin r + \epsilon q(r),
\end{equation}
where $\epsilon$ is an infinitesimal deformation parameter and $p(r)$ and $q(r)$ are smooth functions.
Then the differential equations \eq{sqs4-deq} at the first order of $\epsilon$ reduce to
\begin{eqnarray} \label{inf-deq}
&& p'' - q = (p' + q')\cot r  + \frac{3}{\sin^2 r} (p-q), \xx
&& q'' + q - 2p = 2  p' \cot r - \frac{6}{\sin^2 r} (p-q).
\end{eqnarray}
Define $h \equiv p-q$. Then, $h'' + h' \cot r + \left( 2 - \frac{9}{\sin^2 r} \right) h = 0$.
A general solution is
\begin{equation}\label{gen-sold}
  h(r) = \frac{1}{\sin^3 r} \big( a_1 \cos r + a_2 (4 \sin^2 r + \sin^4 r - 8 ) \big).
\end{equation}
Since $h(r)$ must be a smooth and finite function on $r \in [0, \pi]$, we require $a_1 = a_2 = 0$, i.e., $p=q$.
With this condition, a general solution of Eq. \eq{inf-deq} is
\begin{equation}\label{gen-solh}
  p(r) = q(r) = c_1 \cos r + c_2 (\sin r - r \cos r),
\end{equation}
where $c_1$ and $c_2$ are arbitrary constants.

The existence of the solution \eq{gen-solh} appears as if a new Einstein manifold has been obtained.
At first sight, it seems to contradict the rigidity of the Einstein structure \cite{besse,daga-yang}.
However, it turns out to be a trivial solution, at most equivalent to changing only the size of the four-sphere.
First let us examine the Euler characteristic $\chi(M)$ and the Hirzebruch signature $\tau(M)$ by expanding them
for the infinitesimal deformation \eq{inf-def} up to the first order such that $\chi(M) = 2  + \epsilon \chi_1 (M) + \cdots$ and
$\tau(M) = 0  + \epsilon \tau_1 (M) + \cdots$.
It is easy to see that $\tau_1 (M) = 0$ since $f(r)=g(r)$
in \eq{inf-def}, so $F^{(\pm)i} \wedge F^{(\pm)i} = \pm \frac{3}{2} \Big( \frac{f''}{f} \Big)^2 d \mu$.
It is also easy to verify that $\chi_1 (M) = 0$. Thus the topological numbers do not change under
the deformation \eq{inf-def}, as they should be.

The $SU(2)$ connections \eq{sqs4-gauge} under the infinitesimal deformation \eq{inf-def} are changed as
\begin{eqnarray} \label{def-su2conn}
 A'^{(\pm)1} &=& \frac{1}{4} (-1 \pm \cos r) d \theta \pm \frac{\epsilon}{4} s(r) d \theta, \xx
 A'^{(\pm)2} &=& \frac{1}{4} (-1 \pm \cos r) \sin \theta d \phi \pm \frac{\epsilon}{4} s(r) \sin \theta d \phi, \\
 A'^{(\pm)3} &=& - \frac{1}{2} \cos \theta d \phi + \frac{1}{4} (1 \pm \cos r) (d \psi + \cos \theta d \phi)
 \pm \frac{\epsilon}{4} s(r) (d \psi + \cos \theta d \phi), \nonumber
\end{eqnarray}
where $s(r) = (c_2 r - c_1) \sin r$. Since the added terms linear in $\epsilon$ simply correspond to a replacement
$f(r) = g(r) = \sin r \to \tilde{f}(r)= \tilde{g}(r) = \sin r  + \epsilon p(r)$ in the metric \eq{sq-s4},
we can use Eq. \eq{sqs4-f} to calculate the curvatures for the deformed connections \eq{def-su2conn}.
The result is\footnote{The reason why $\chi_1 (M) = 0$ is that the volume form also has a linear correction,
$d\mu = d\mu (\mathbb{S}^4) \left (1 + 3 \epsilon \frac{p(r)}{\sin r} \right)$, which exactly cancels
the linear term in \eq{def-scur}.}
\begin{equation}\label{def-scur}
  F'^{(\pm)i} = \left(\frac{1}{2} - \epsilon c_2 \right) \zeta^i_\pm.
\end{equation}

The result \eq{def-scur} shows that the smooth deformations \eq{inf-def} preserve the Einstein structure since
$f^{ij}_{(\pm\mp)}$ identically vanish.
It turns out that even the Einstein-Weyl structure \cite{pps-im1996,taka-tjm,bonn-cqg}
has not been deformed by the solutions \eq{gen-solh}.
The Einstein-Weyl structure deformations are captured by
${\widetilde f}^{ij}_{(\pm\pm)} \equiv f^{ij}_{(\pm\pm)}
- \frac{1}{3} \delta^{ij} \big(  f^{kl}_{(\pm\pm)} \delta^{kl} \big)$ in Eq. \eq{weyl-riem}.
For the curvature \eq{def-scur}, $f^{ij}_{(\pm\pm)} = \left(\frac{1}{2} - \epsilon c_2 \right) \delta^{ij}$
and so ${\widetilde f}^{ij}_{(\pm\pm)} = 0$. This means that the Einstein metric determined by the functions
in \eq{inf-def} is still conformally flat. This in turn implies that the Einstein-Weyl structure of a round four-sphere
is preserved, while maintaining the Einstein structure.
It is allowed to change at most only the size of the sphere.
Such a constant scaling is regarded as trivial (see Theorem $\mathbf{5.41}$ in \cite{besse}).
This property was crucial in \cite{hkky,kky} to realize the instanton-induced inflation of four-dimensional spacetime
and the dynamical compactification of internal space in higher-dimensional Einstein-Yang-Mills theory.

\section{Taub-NUT/Bolt Instantons}

Our formalism proves particularly useful in the Taub-NUT and Taub-Bolt spaces.
The Taub-NUT solution of the vacuum Einstein equations exhibits intriguing properties \cite{taub,nut,misner}.
The metric for the general Taub-NUT space with a Lorentzian signature is given by
\be \la{tnut-l}
ds^2 = - \widetilde{f}(r) (dt - 2n \cos \theta d\phi)^2 + \frac{1}{\widetilde{f}(r)} dr^2 + (r^2 + n^2) (d\theta^2 + \sin^2 \theta d\phi^2),
\ee
where
\be \la{l-pot}
\widetilde{f}(r) = \frac{r^2 - n^2 - 2mr}{r^2 + n^2}.
\ee
Here $m$ is the ADM (gravitoelectric) mass of the solution, and $n$ is the NUT charge, regarded as a gravitomagnetic mass.
Therefore the Taub-NUT solution can be considered to be a gravitational dyon.
Note that the function \eq{l-pot} can be written in an intriguing form
\be \la{sep-pot}
\widetilde{f}(r) = 1 - \frac{m+in}{r+in} - \frac{m-in}{r-in}.
\ee
Its physical interpretation becomes more clear in Euclidean formulation.

The general Euclidean Taub-NUT family is obtained by performing the Wick rotation which transforms both time and the NUT charge:
\be \la{wick-rot}
t \to i \tau, \quad n \to i n.
\ee
The resulting metric is known as the Taub-Bolt space \cite{gh-cmp1979,tbolt-page,gibb-pope} and takes the form
\be \la{tnut-e}
ds^2 = f(r) (d\tau - 2n \cos \theta d\phi)^2 + \frac{1}{f(r)} dr^2 + \rho^2 (d\theta^2 + \sin^2 \theta d\phi^2),
\ee
where
\be \la{e-pot}
f(r) = \frac{r^2 + n^2 - 2mr}{\rho^2} = 1 - \frac{m+n}{r+n} - \frac{m-n}{r-n}, \quad \rho^2 = r^2 - n^2 \;\; (r \geq |n|).
\ee
Note that setting $n = \pm m$ changes the function \eq{e-pot} to $f(r) = \frac{r-m}{r+m}$
and leads to self-dual Taub-NUT metrics \cite{gh-cmp1979,egh-report}
\be \la{sd-tnut}
ds^2 = \left( \frac{r-m}{r+m} \right) (d\tau \mp 2m \cos \theta d\phi)^2 + \left( \frac{r+m}{r-m} \right) dr^2 + (r^2 - m^2) (d\theta^2 + \sin^2 \theta d\phi^2).
\ee
One may further shift the coordinate $r-m \to r $ such that $r \geq 0$ and $f(r) = \frac{r-m}{r+m} \to V(r)^{-1} =  \left( 1 + \frac{2m}{r} \right)^{-1}$ to get the metric of Gibbons-Hawking gravitational instantons \cite{hawk-sdi,gibb-hawk-gi}:
\be \la{gh-gi}
ds^2 = V(r)^{-1} (d\tau \mp a)^2 +  V(r) \big( dr^2 + r^2 (d\theta^2 + \sin^2 \theta d\phi^2) \big),
\ee
where $ V(r) = 1 + \frac{2m}{r}$ and the gauge field, $a = - 2n \cos \theta d\phi$, is clearly
that of a monopole \cite{gross-perry}
\be \la{pot-mono}
\mathbf{B} = \mathbf{\nabla} \times \mathbf{a} = \frac{2n}{r^2} \hat{\mathbf{r}}.
\ee

The Euclidean vierbein for the metric \eq{tnut-e} is
\be \la{vierbein}
e^1 = \frac{dr}{\sqrt{f(r)}}, \quad e^2 = \rho d\theta, \quad  e^3 = \rho \sin \theta d\phi,
\quad e^4 = \sqrt{f(r)} (d\tau - 2n \cos \theta d\phi).
\ee
After getting spin connections from the orthnormal basis in \eq{vierbein}, one can identify the self-dual
and anti-self-dual $SU(2)$ gauge fields in \eq{ga-split}:
\bea \la{tnut-conn}
A^{(\pm)1} &=& \frac{1}{2} \left( \omega_{23} \pm \omega_{14} \right) = - \frac{1}{2} \left( \frac{1}{\sqrt{\rho}} \cot \theta e^3 + \frac{n \sqrt{f}}{\rho^2} e^4
\pm  \frac{1}{\sqrt{f}} \left(  \frac{m}{\rho^2} - \frac{2 n^2}{\rho^4} (r-m) \right) e^4 \right), \xx
A^{(\pm)2} &=& \frac{1}{2} \left( \omega_{31} \pm \omega_{24} \right) = \frac{\sqrt{f}}{2\rho^2} (r \mp n) e^3, \\
A^{(\pm)3} &=& \frac{1}{2} \left( \omega_{12} \pm \omega_{34} \right) =  - \frac{\sqrt{f}}{2\rho^2} (r \mp n) e^2. \nonumber
\eea
One can determine the $SU(2)$ curvatures defined by \eq{sd-2f}
\bea \la{tnut-curv}
F^{(\pm)1} &=& \frac{m \pm n}{(r \pm n)^3} ( e^2 \wedge e^3  \pm e^1 \wedge e^4), \xx
F^{(\pm)2} &=& \frac{m \pm n}{2 (r \pm n)^3} ( e^1 \wedge e^3  \pm e^4 \wedge e^2), \\
F^{(\pm)3} &=& \frac{m \pm n}{2 (r \pm n)^3} ( e^1 \wedge e^2  \pm e^3 \wedge e^4). \nonumber
\eea
It is obvious that $F^{(+)i}$ are self-dual two-forms and describe a gravitational instanton while
$F^{(-)i}$ are anti-self-dual two-forms and describe a gravitational anti-instanton.
It is interesting to notice that gravitational instanton and anti-instanton in the Taub-Bolt space \eq{tnut-e}
are simply related to each other by a sign flip of the nut charge:
\be \la{+--corr}
F^{(-)i}_{ab}(n) = F^{(+)i}_{ab} ( n \to -n).
\ee
If one sets $n=m$ (a nut), $F^{(-)i}$ identically vanish, so only self-dual instanton exists.
But, if one sets $n = - m$ (an anti-nut), $F^{(+)i}$ instead dissapear, so only anti-self-dual instanton survives.
This feature may explain why the function $f(r)$ in \eq{e-pot} contains two kinds of contributions.

The instanton densities \eq{density} for the self-dual ($n=m$) and anti-self-dual ($n=-m$) cases are given by
\be \la{ins-den-tnut}
F^{(\pm)i} \wedge F^{(\pm)i} = \pm 12 m^2 \frac{1}{(r+m)^6} d \mu,
\ee
where $d\mu = (r^2 - m^2) \sin \theta dr d\theta d\phi d\tau$ is the volume form of the Taub-NUT space.
Then the instanton number for the density \eq{ins-den-tnut} can easily be calculated to be
\be \la{ins-num-tnut}
I^{(\pm)} \equiv \frac{1}{4\pi^2} \int_M F^{(\pm)i} \wedge F^{(\pm)i} = \pm \frac{12 m^2 \gamma}{\pi} \int_m^\infty \frac{r^2 - m^2}{(r+m)^6} d r = \pm 1,
\ee
where we used the fact \cite{gh-cmp1979,gibb-pope} that the period of the coordinate $\tau \in \mathbb{S}^1$ is $\gamma = \int_0^\gamma d\tau = 8 \pi m$.

For the general Taub-Bolt metric \eq{tnut-e} which is not self-dual or anti-self-dual, the radial coordinate $r$ must range
from a zero in $f(r)$ to infinity, but it must avoid the values $\pm n$, which are both curvature singularities.
It turns out \cite{tbolt-page} that one must take $m = \frac{5}{4}n$ and the ranges
\be \la{coord-range}
2n \leq r < \infty, \quad 0 \leq \theta \leq \pi, \quad 0 \leq \phi \leq 2 \pi, \quad 0 \leq \tau \leq 8 \pi n.
\ee
Then the Euler characteristic $\chi(M)$ is
\bea \la{euler-tbolt}
\chi(M) &=& I^{(+)} - I^{(-)} \xx
&=& 24 n \left( (m+n)^2  \int_{2m}^\infty \frac{r^2 - m^2}{(r+m)^6} d r + (m-n)^2  \int_{2m}^\infty \frac{r^2 - m^2}{(r-m)^6} d r \right) \xx
&=& 24 n \left( \frac{1}{32 n} + \frac{5}{96 n} \right) =  2,
\eea
and the Hirzebruch signature $\tau(M)$ is given by
\bea \la{sign-tbolt}
\tau(M) &=& \frac{2}{3} (I^{(+)} + I^{(-)}) - \frac{2}{3k} + \frac{(k-1)(k-2)}{3k} \xx
&=&  16 n \left( \frac{1}{32 n} - \frac{5}{96 n} \right) - \frac{2}{3} = -1.
\eea
The extra terms in $\tau(M)$ are the contributions from a boundary $\partial M = \mathbb{S}^3/\mathbb{Z}_k$ (a cyclic lens space)
and $k=1$ for the Taub-Bolt space,
so the $\eta$-invariant $\eta_S (\partial M)=\frac{(k-1)(k-2)}{3k}$ identically vanishes \cite{gpr-index,gp-index}.

The Hitchin-Thorpe inequality \eq{hitchin-thorpe} can be extended to asymptotically locally flat (ALF) spaceces \cite{gp-index} such as the Taub-NUT/Bolt spaces \eq{tnut-e} using Eqs. \eq{euler-tbolt} and \eq{sign-tbolt}
\begin{equation}\label{hitchin-thorpe}
  \chi (M) \geq \frac{3}{2} \left| \tau (M) + 1 - \frac{k}{3} \right|,
\end{equation}
where the equality is obtained in the case of half-flat spaces, i.e. either $I^{(+)} = 0$ or $I^{(-)} = 0$.
Since both Taub-NUT and Taub-Bolt spaces are Ricci-flat and satisfy the same ALF condition,
one might wonder whether there could be a smooth deformation connecting the two spaces.
However, the Hitchin-Thorpe inequality implies that one cannot find a small perturbation of a solution
to give a nearby solution with a different topology \cite{rev-anderson,gp-index}.
In other words, the Taub-NUT space $(\chi=1, \; \tau = 0)$ carries
an Einstein structure different from the Taub-Bolt space $(\chi=2, \; \tau = \pm 1)$
even though they share the same ALF boundary condition.
The same conclusion may also be applied to the Euclidean Schwarzschild $(\chi=2, \; \tau = 0)$ and
Kerr $(\chi=2, \; \tau = 0)$ solutions. In contrast to the previous example, in this case,
their asymptotic boundary conditions are different although they are both Ricci-flat and share the same topological numbers.
This issue will be addressed elsewhere.

\section{Moduli Space of Einstein Metrics}

The moduli space $\mathcal{E}(M)$ of Einstein structures on a compact four-manifold $M$ is the quotient \cite{besse,rev-anderson}
\be \la{emoduli-space}
\mathcal{E}(M) = \{\mathrm{Einstein \; metrics \; on \; \textit{M}} \}/(\mathrm{Diff}(M) \times  \mathbb{R}^+),
\ee
where $\mathrm{Diff}(M)$ is diffeomorphisms of $M$ and the positive real $\mathbb{R}^+$ acts by rescaling.
The (local) triviality of this space is essentially originated from the fact that
arbitrary deformations of the Ricci tensor satisfy the so-called Palatini identity:
\be \la{palatini}
\delta R_{\mu\nu} = \nabla_\rho \delta {\Gamma^\rho}_{\mu\nu} - \nabla_\nu \delta {\Gamma^\rho}_{\mu \rho},
\ee
where $\delta {\Gamma^\rho}_{\mu\nu}$ denotes an arbitrary metric deformation of Christoffel symbols and $\nabla_\rho$ indicates
a covariant differentiation defined by an undeformed metric.
If one contracts the (undeformed) inverse metric density, $\sqrt{g} g^{\mu\nu}$, with Eq. \eq{palatini},
the result becomes a total derivative,
\be \la{palatini-total}
\sqrt{g} g^{\mu\nu} \delta R_{\mu\nu} = \partial_\mu ( \sqrt{g} A^\mu ),
\ee
where $A^\mu = g^{\rho \sigma} \delta {\Gamma^\mu}_{\rho\sigma} - g^{\rho \mu} \delta {\Gamma^\sigma}_{\rho\sigma}$.
The identity \eq{palatini-total} leads to an important result \cite{besse,rev-anderson} that,
for a small deformation $g \to g + \delta g$,
an Einstein metric $g$ is critical for the functional $\mathcal{F}(g)$ defined by
\be \la{eins-funal}
 \mathcal{F}(g) = vol(g)^{-1/2} \int_M R(g) d\mu,
\ee
where $vol(g) = \int_M d\mu$ and $d\mu = \sqrt{g} d^4 x$ is the volume form introduced in \eq{intersection}.
This property is absent in gauge theories.
The same identity also holds for the Lie derivative:
\be \la{palatini-lie}
\mathcal{L}_X R_{\mu\nu} = \nabla_\rho \left( \mathcal{L}_X {\Gamma^\rho}_{\mu\nu} \right) - \nabla_\nu  \left(\mathcal{L}_X {\Gamma^\rho}_{\mu \rho} \right),
\ee
where $X = X^\mu \partial_\mu$ denotes any vector field on a manifold $M$. This implies that some deformations in \eq{palatini}
are actually the diffeomorphisms generated by the vector field $X$, so will be trivialized in the moduli space $\mathcal{E}(M)$.

Using Eqs. \eq{coeff-1} and \eq{coeff-2}, the functional \eq{eins-funal} in our approach can be written as
\be \la{funcal-our}
\mathcal{F}(A^{(\pm)}, \zeta_\pm) = 2 \int_M \left( F^{(+)i} \wedge \zeta^i_+ - F^{(-)i} \wedge \zeta^i_- \right).
\ee
Since the functional \eq{eins-funal} is invariant under scaling, we have assumed all metrics are normalized
to have unit volume, i.e.
\be \la{normal-vol}
vol(g) = \frac{1}{2} \int_M  \zeta^i_+ \wedge \zeta^j_+  = - \frac{1}{2} \int_M  \zeta^i_- \wedge \zeta^j_- = \delta^{ij}.
\ee
Since we apply the tetrad formalism where the general coordinate invariance is manifest,
the corresponding local symmetry in \eq{emoduli-space} is replaced by the local Lorentz symmetry in \eq{spin4}
and the scaling symmetry $\mathbb{R}^+$ is now fixed by the condition \eq{normal-vol}.
Including the constraints \eq{normal-vol}, the functional \eq{funcal-our} is modified as
\be \la{funcal-mod}
\mathcal{F}_\alpha (A^{(\pm)}, \zeta_\pm) = \int_M \left[ 2 \left( F^{(+)i} \wedge \zeta^i_+ - F^{(-)i} \wedge \zeta^i_- \right)
- \alpha^{ij}_+ (\zeta^i_+ \wedge \zeta^j_+ - 2 d\mu \delta^{ij})
+ \alpha^{ij}_- (\zeta^i_- \wedge \zeta^j_- + 2 d\mu \delta^{ij})  \right],
\ee
where the Lagrange multipliers $\alpha^{ij}_\pm$ are symmetric quantities.
The variation of $\mathcal{F}_\alpha$ gives the following equations of motion:
\bea \la{eom1-fa}
&&   F^{(+)i}  - \alpha^{ij}_+ \zeta^j_+  = 0,
\qquad F^{(-)i}  - \alpha^{ij}_- \zeta^j_-  = 0, \\
\la{eom2-fa}
&& D^{(+)} \zeta^i_+ = D^{(-)} \zeta^i_- = 0, \\
\la{eom3-fa}
&& \zeta^i_+ \wedge \zeta^j_+  =  - \zeta^i_- \wedge \zeta^j_- = 2 d\mu \delta^{ij},
\eea
where $D^{(\pm)} \zeta^i_\pm = d \zeta^i_\pm   - 2 \varepsilon^{ijk} A^{(\pm)j} \wedge \zeta^k_\pm$.
We see that Eq. \eq{eom1-fa} reproduces \eq{sde} with $\alpha^{ij}_\pm = f^{ij}_{(\pm\pm)}$.
Note also that $F^{(\pm)i} \wedge  \zeta^i_\pm  = \pm 2 \alpha^{ij}_\pm \delta^{ij} d\mu$,
so $\lambda = 2 \alpha^{ij}_\pm \delta^{ij}$ is the cosmological constant of the Einstein manifold $M$.
Eq. \eq{eom2-fa} corresponds to the torsion-free condition \eq{t-free} and it can be solved
algebraically for $ A^{(\pm)i}$.  Therefore an Einstein manifold with or without a cosmological constant
is critical for the functional $\mathcal{F}_\alpha $.

We now study the equations governing perturbations of solutions to Eqs. \eq{eom1-fa}-\eq{eom3-fa}.
Let us denote the linear perturbations by $B^{(\pm)i}$ and $\xi^i_\pm$ and make the replacement
\be \la{perturb}
A^{(\pm)i} \to A'^{(\pm)i} = A^{(\pm)i} + B^{(\pm)i}, \qquad
\zeta^i_\pm \to \zeta'^i_\pm = \zeta^i_\pm + \chi^i_\pm.
\ee
If the functional $\mathcal{F}_\alpha$ is expanded to second-order in the perturbations,
one can find the following result
\bea \la{pert-exp}
\mathcal{F}_\alpha (A'^{(\pm)}, \zeta'_\pm) &\approx & \mathcal{F}_\alpha (A^{(\pm)}, \zeta_\pm) \xx
&& + \int_M \left[ 2 \left( D^{(+)} B^{(+)i} \wedge \chi^i_+ - D^{(-)} B^{(-)i} \wedge \chi^i_- \right) \right. \xx
&& - 2 \varepsilon^{ijk} \left( \zeta^i_+ \wedge B^{(+)j} \wedge B^{(+)k} - \zeta^i_- \wedge B^{(-)j}
\wedge B^{(-)k}    \right) \\
&& \left.
- \alpha^{ij}_+ \chi^i_+ \wedge \chi^j_+  + \alpha^{ij}_- \chi^i_- \wedge \chi^j_-
- 2 \left( \delta \alpha^{ij}_+ \zeta^i_+ \wedge \chi^j_+  - \delta \alpha^{ij}_- \zeta^i_- \wedge \chi^j_- \right)
\right], \nonumber
\eea
where the equations of motion, \eq{eom1-fa}-\eq{eom3-fa}, are used for the first-order terms.
We observe that the Lagrange multiplier terms can be consistent only by choosing the ansatz $\chi^i_\pm = \zeta^i_\mp$.
Using \eq{intersection}, the second-order perturbations reduce to
\bea \la{2nd-order}
\mathcal{F}^{(2)}_\alpha &=& \int_M \left[ 2 \left( \zeta^i_- \wedge D^{(+)} B^{(+)i}
- \zeta^i_+ \wedge D^{(-)} B^{(-)i} \right) \right. \xx
&& \left. - 2 \varepsilon^{ijk} \left( \zeta^i_+ \wedge B^{(+)j} \wedge B^{(+)k} - \zeta^i_- \wedge B^{(-)j}
\wedge B^{(-)k} \right) + 2 \lambda d\mu  \right].
\eea
One can show that local $SU(2)_\pm$ transformations $B^{(\pm)i} =  D^{(\pm)} \phi^i_\pm$ trivialize
the second-order functional $\mathcal{F}^{(2)}_\alpha$ if the gauge parameters satisfy
$\phi_\pm^i \phi^j_\pm = \phi_\pm^2 \delta^{ij}$ and $ \phi_+^2 + \phi_-^2
= \frac{1}{4}$.

The above test of the gauge solution gives the prediction that the ansatz $\chi^i_\pm = \zeta^i_\mp$
is also required in general cases. It implies that, in contrast to gauge theories,
perturbations in gravitational theory require an opposite self-duality
compared to the self-duality of the unperturbed curvature.
A similar behavior was also observed in \cite{c-torre} (see the last paragraph in Sect. IV.).
The equations of motion derived from the second-order functional \eq{2nd-order} are
\be \la{2nd-eom}
 D^{(\pm)} \zeta^{i}_\mp - 2 \varepsilon^{ijk}  B^{(\pm)j} \wedge \zeta^{k}_\pm = 0.
\ee
If a parent Einstein manifold is given, the connections $D^{(\pm)}$ as well as
the self-dual 2-form bases $\zeta^{i}_\pm$ are determined. Therefore one can solve Eq. \eq{2nd-eom} algebraically
for $B^{(\pm)i}$. A trivial solution from the gauge transformation is
$\zeta^{i}_\mp - 2 \varepsilon^{ijk}  \phi^j_\pm \zeta^{k}_\pm = \eta^i_\pm$ where
$\eta^i_\pm$ are covariantly constants.
Our interest is whether there is a nontrivial solution to \eq{2nd-eom} other than the gauge transformation.
Currently, we do not have a definite answer. However, our attempt using some examples leads to the conclusion that
Eq. \eq{2nd-eom} has a structure (mixing opposite self-dualities) that is difficult to admit solutions.
Supported on the putative rigidity for the moduli space of Einstein metrics too, we conjecture
that there are no nontrivial solutions to Eq. \eq{2nd-eom} and the moduli space of gauge-inequivalent solutions
is discrete, i.e., zero-dimensional. We hope to clarify this issue in the near future.

\section{Discussion}

We expect that $\mathbb{C}\mathbb{P}^2$ will exhibit a similar feature as $\mathbb{S}^4$
for the Einstein deformations \cite{besse,daga-yang}.
The Theorem $\mathbf{2.2}$ also implies that any continuous deformations of $\mathbb{C}\mathbb{P}^2$
will break the Einstein structure.
However, some details are worth studying for the case of $\mathbb{C}\mathbb{P}^2$.
Interestingly, $\mathbb{S}^2 \times \mathbb{S}^2$ would reveal a very different aspect of Einstein structures
since all two-dimensional Riemannian manifolds are Einstein spaces. So any smooth deformations of $\mathbb{S}^2$ (e.g., ellipsoids)
will not change the underlying Einstein structure. This reasoning is consistent with the Proposition 12.75
(the case $SU(m)/SO(m)$ for $m=2$) in \cite{besse}.

Let us consider from a physical perspective why the deformation of Einstein structures is so rigid.
The Einstein equations with matters take a more general form
\begin{equation}\label{matter-eins-eq}
  R_{\mu\nu} - \frac{1}{2} g_{\mu\nu} R + \lambda g_{\mu\nu} = \frac{8 \pi G}{c^4} T_{\mu\nu},
\end{equation}
where the right-hand side refers to the energy-momentum tensor of matters.
Without matters, the Einstein equations reduce to $R_{\mu\nu} = \lambda g_{\mu\nu}$,
whose solutions are called Einstein manifolds. Therefore, the moduli space of Einstein structures refers to
the space of all possible solutions of Einstein equations without matters (i.e., in vacuum).
This system is not elliptic due to the diffeomorphism invariance of the equations.
However, for a suitable gauge choice (e.g., the de Donder gauge), the restricted equations become elliptic \cite{rev-anderson}.
To uniquely determine a solution to Einstein equations, it is necessary to specify suitable initial conditions
or boundary conditions. The uniqueness of solutions can depend on several factors such as the choice of boundary conditions
and the presence of symmetries, which are usually asymptotic and global conditions.
Given a solution under the specified condition, local deformations of the given solution cannot change
such an asymptotic and/or global condition underling the solution.
Therefore, the existence of a nontrivial moduli space of Einstein structures depends on whether a new Einstein metric
can be generated by local deformations while maintaining such asymptotic and global conditions.
Considering that spacetime behaves like a metrical elasticity  with tension, delicate balances are required to maintain
such deformations. For example, the Page metric is obtained as a limiting case of
the Kerr-de Sitter solution \cite{page,gh-cmp1979}.
This may be a reason why the Page space could be an Einstein metric even though it is nonhomogeneous.
Furthermore, it is a nontrivial $\mathbb{S}^2$ bundle over $\mathbb{S}^2$ which is more flexible against (local) deformations
than other Einstein spaces. For this reason, it seems very unlikely that one parameter family of Einstein structure
deformations exists.

Recent progresses in four-dimensional topology has shown that many four-dimensional manifolds have exotic differential
structures \cite{dona-kron}. However, the Poincar\'e conjecture in dimension 4 remains to be an open problem,
that is, whether there exists exotic differential structures on the four-sphere.
According to the Freedman’s classification \cite{freedman}, two smooth, compact, simply-connected four-manifolds
are homeomorphic if and only if they are both spin or both non-spin and their Euler characteristics and signatures are equal.
Four-sphere is consisted of an $SU(2)_+$ instanton and an $SU(2)_-$ anti-instanton which determines $\chi = 2$ and $\tau = 0$
\cite{oh-yang,hkky} as shown in Fig. \ref{top-num}. Thus it is impossible to change the instanton composition
in order for two four-spheres to have the same smooth structure.
Perhaps only deforming the arrangement of two instantons is allowed. But it corresponds to the deformation of
Einstein structures on $\mathbb{S}^4$. Therefore the exotic differential structures on the four-sphere may be deeply related
to its Einstein structures \cite{kots-1998}.

\section*{Acknowledgments}

We are deeply grateful to Prof. Youngjoo Chung for his help on the solution of differential equations.
HSY thanks Jongmin Park for his help on some calculations in Section 5.
This work was supported by the National Research Foundation of Korea (NRF) with grant numbers
NRF-2019R1A2C1007396 (KK) and NRF-2018R1D1A1B0705011314 (HSY).
JH was supported by the Ministry of Education through the Center for Quantum Spacetime
(CQUeST) of Sogang University (NRF-2020R1A6A1A03047877) and NRF-2020R1A2C1014371.
We acknowledge the hospitality at APCTP where part of this work was done.

\appendix

\section{Strong version of the Hitchin-Thorpe inequality}

The Hitchin–Thorpe inequality \cite{hitchin1974,besse} provides an important obstruction to the existence
of Einstein metrics on four-manifolds. In this appendix, we will prove the strong version of the Hitchin-Thorpe
inequality \eq{stineq-hitchin}. Here we assume that a compact Einstein manifold has nonnegative sectional curvature.

The curvature tensors for an Einstein manifold are given by \eq{em-riem} where
the coefficients $f^{ij}_{(\pm\pm)}$ are real symmetric $3 \times 3$ matrices.
Therefore they can be diagonalized by orthogonal transformations as \cite{chy-grg}\footnote{Since $SU(2)_+$
and $SU(2)_-$ parts are independent of each other, they can be diagonalized separately.}
\begin{equation}\label{2-matrices}
  (A_\pm)_{ij} \equiv f^{ij}_{(\pm\pm)} = \mathrm{diag} (a_\pm^1, a_\pm^2, a_\pm^3).
\end{equation}
The diagonal elements satisfy the relation
\begin{equation}\label{pseudo0}
a_+^1 + a_+^2 + a_+^3= a_-^1 + a_-^2 + a_-^3 = \frac{R}{8}
\end{equation}
according to Eqs. \eq{1-more} and \eq{ricci}.

We regard the curvature tensor $R$ as a symmetric linear transformation of the bundle $\Omega^2 (M)$ of 2-forms
defined by \cite{hitchin1974}
\begin{equation}\label{curv-map}
  R(e^a \wedge e^b) = R_{ab} = \frac{1}{2} R_{abcd} e^c \wedge e^d
\end{equation}
relative to a local orthonormal basis $\{e^a \}$ of the 1-forms.
We choose the basis of $\Omega^2 (M)$ as $\{e^1 \wedge e^2, e^3 \wedge e^1, e^2 \wedge e^3, e^3 \wedge e^4, e^2 \wedge e^4,
e^1 \wedge e^4\}$ according to the canonical pairing \eq{tri-2form}. Then $R$ takes the form \cite{hitchin1974}
\begin{equation}\label{r-matrix}
  R = \left(
        \begin{array}{cc}
          A & B \\
          B & A \\
        \end{array}
      \right)
\end{equation}
where
\begin{equation}\label{ab-matrix}
  A = \left(
        \begin{array}{ccc}
          \lambda_1 & 0 & 0 \\
          0 & \lambda_2 & 0 \\
          0 & 0 & \lambda_3 \\
        \end{array}
      \right), \qquad
  B = \left(
        \begin{array}{ccc}
          \mu_1 & 0 & 0 \\
          0 & \mu_2 & 0 \\
          0 & 0 & \mu_3 \\
        \end{array}
      \right),
\end{equation}
and $\lambda_1 = a_+^3 + a_-^3 , \; \lambda_2 = a_+^2 + a_-^2,
\; \lambda_3 = a_+^1 + a_-^1$, and  $\mu_1 = a_+^3 - a_-^3 , \; \mu_2 = a_+^2 - a_-^2,
\; \mu_3 = a_+^1 - a_-^1$.
The condition for the nonnegative sectional curvature means that $\lambda_i \geq 0, \, i=1,2,3$.
Eq. \eq{pseudo0} imposes the conditions $\lambda_1 + \lambda_2 + \lambda_3 = \frac{R}{4} > 0$
and $\mu_1 + \mu_2 + \mu_3  = 0$. Thus we set $\mu_3 = - (\mu_1 + \mu_2)$ and
$\lambda_1 \geq \lambda_2 \geq \lambda_3 \geq 0$.
We regard $\vec{\lambda} =( \lambda_1, \lambda_2, \lambda_3)$ and
$\vec{\mu} =( \mu_1, \mu_2, \mu_3)$ as vectors in $\mathbb{R}^3$.
Since the vector $\vec{\mu} =( \mu_1, \mu_2, \mu_3)$ is constrained by $\mu_1 + \mu_2 + \mu_3  = 0$,
if $\vec{\lambda}, \; \vec{\mu}$ are nonzero, then the angle $\theta$ between them must satisfy
\begin{eqnarray} \label{angle-ineq}
  \cos \theta &=& \frac{\vec{\lambda} \cdot \vec{\mu}}{|\vec{\lambda}| |\vec{\mu}|}
  \leq \frac{1}{|\vec{\lambda}| \sqrt{6 \mu_1 \mu_2}}
  \big( (\lambda_1 - \lambda_3) \mu_1 + (\lambda_2 - \lambda_3) \mu_2 \big)
  \leq \sqrt{\frac{2}{3}} \frac{\sqrt{(\lambda_3 - \lambda_1) (\lambda_3 - \lambda_2)}}{|\vec{\lambda}|} \xx
  &\leq & \sqrt{\frac{2}{3}},
\end{eqnarray}
using the inequality between arithmetic and geometric means.
Therefore
$$ \vec{\lambda} \cdot \vec{\mu} \leq \sqrt{\frac{2}{3}} |\vec{\lambda}| |\vec{\mu}|.$$
From the inequality $|\vec{\lambda}|^2 +  |\vec{\mu}|^2 \geq 2 |\vec{\lambda}| |\vec{\mu}|$, we get
\begin{equation} \label{cru-ineq}
  |\vec{\lambda}|^2 +  |\vec{\mu}|^2 \geq 2 \sqrt{\frac{3}{2}} (\vec{\lambda} \cdot \vec{\mu}).
\end{equation}

Note that
$$ |\vec{\lambda}|^2 +  |\vec{\mu}|^2 = 4 \pi^2 \rho_\chi, \qquad
(\vec{\lambda} \cdot \vec{\mu}) = 3 \pi^2 \rho_\tau,$$
where the densities were defined in \eq{def-de} and \eq{def-dh}.
Then the inequality \eq{cru-ineq} can be written as
\begin{equation}\label{sht-ineq}
  \rho_\tau \leq \left( \frac{2}{3} \right)^{\frac{3}{2}} \rho_\chi.
\end{equation}
Integrating and using the expressions for the Euler characteristic and Hirzebruch signature
in \eq{den-euler} and \eq{den-hirze}, we obtain the inequality \eq{stineq-hitchin}.
Incidentally, the Hitchin-Thorpe inequality \eq{hitchin-thorpe} can also be derived
from the inequality
$$4 \pi^2 \left(\rho_\chi \pm \frac{3}{2}  \rho_\tau \right) = (\vec{\lambda} \pm  \vec{\mu})^2 \geq 0.$$
The case of nonpositive sectional curvature is similar.

The Hitchin-Thorpe inequality can be further refined by rewriting the topological invariants
in Eqs. \eq{dec-euler} and \eq{dec-hirze} as follows:
\begin{eqnarray} \label{ref-euler}
\chi(M) &=& \frac{1}{2 \pi^2} \int_M \left( (\widetilde{f}^{ij}_{(++)})^2 + (\widetilde{f}^{ij}_{(--)})^2
+ \frac{R^2}{96} \right) d\mu, \\
\label{ref-hirze}
\tau(M) &=& \frac{1}{3 \pi^2} \int_M \left( (\widetilde{f}^{ij}_{(++)})^2 - (\widetilde{f}^{ij}_{(--)})^2 \right) d\mu,
\end{eqnarray}
where $\widetilde{f}^{ij}_{(\pm\pm)} = f^{ij}_{(\pm\pm)} - \frac{R}{24} \delta^{ij}$ are the coefficients
in the expansion of Weyl tensors in Eq. \eq{weyl-riem}. Then we get
\begin{equation}\label{ref-hittho}
\chi(M) \pm \frac{3}{2}  \tau(M) = \frac{1}{\pi^2} \int_M \left( \widetilde{f}^{ij}_{(\pm \pm)} \right)^2 d\mu
+ \frac{R^2}{192 \pi^2} \mathrm{vol}(M),
\end{equation}
where $\mathrm{vol}(M)$ is the volume of an Einstein manifold $M$.
Since the scalar curvature is given by $R = 4 \lambda$, Eq. \eq{ref-hittho} can be written as the form
\begin{equation}\label{ref-hittho2}
\chi(M) - \frac{\lambda^2}{12 \pi^2} \mathrm{vol}(M) \geq \frac{3}{2}  |\tau(M)|,
\end{equation}
where the equality holds only if $\widetilde{f}^{ij}_{(++)} = 0$ or $\widetilde{f}^{ij}_{(--)} = 0$, i.e.,
$M$ is conformally half-flat (such as $\mathbb{S}^4$ and $\mathbb{C}\mathbb{P}^2$).
Indeed, as $\frac{\lambda^2}{12 \pi^2} \mathrm{vol}(M) = 2$ and $\frac{3}{2}$ for $\mathbb{S}^4$ and $\mathbb{C}\mathbb{P}^2$, respectively, the relation \eq{ref-hittho2} holds with the equality.
However, for other spaces, we have, for example,
\begin{equation}\label{2weyl-coeff}
  \widetilde{f}^{ij}_{(\pm\pm)} = \left\{
                                    \begin{array}{ll}
                                      \frac{1}{6 R^2} \, \mathrm{diag} (-1, -1, 2) & \hbox{for $\mathbb{S}^2 \times \mathbb{S}^2$;} \\
                                      \pm \left( \psi^i_{(\pm)} - \frac{1}{3} ( 2 \psi^1_{(\pm)} + \psi^3_{(\pm)} ) \right)
\delta^{ij} & \hbox{for Page},
                                    \end{array}
                                  \right.
\end{equation}
where we used the result \eq{page-f} and $\psi^1_{(\pm)} = \psi^2_{(\pm)}$ for the Page space.
Since both spaces have $\chi = 4$ and $\tau = 0$ as shown in Fig. \ref{top-num},
the inequality \eq{ref-hittho2}, in particular, implies that $\mathrm{vol}(\mathrm{Page}) \leq \frac{48 \pi^2}{\lambda^2}$.
One can see from Eq. \eq{2weyl-coeff} that the (anti-)self-dual Weyl tensors for the Page space have only
two distinct eigenvalues $\left( \pm \frac{1}{3} (\psi^1_{(\pm)} - \psi^3_{(\pm)}),
\; \mp \frac{2}{3} (\psi^1_{(\pm)} - \psi^3_{(\pm)} ) \right)$ at each point.
Therefore, according to Proposition 5 in \cite{v21} (see also Theorem $\mathbf{11.81}$ in \cite{besse}),
the Page metric \eq{page} must be locally conformal to a K\"ahler metric
although it is not a K\"ahler metric in itself.

The inequality \eq{ref-hittho2} plays an important role to find infinitely many compact simply connected
smooth four-manifolds which do not admit Einstein metrics, but nevertheless
satisfy the strict Hitchin-Thorpe inequality $\chi > \frac{3}{2}|\tau|$ \cite{lebrun96}.

\end{document}